\documentclass[reqno]{amsart}
\pdfoutput=1

\usepackage{amsmath,amsthm,amssymb,amsfonts,amscd}
\usepackage[]{hyperref}
\usepackage{booktabs} 
\usepackage{graphicx}
\usepackage{color}
\usepackage{url} 
\usepackage[small,nohug]{diagrams} 
\diagramstyle[labelstyle=\scriptstyle]

\newtheorem{thm}{Theorem}[section]
\newtheorem{prop}[thm]{Proposition}
\newtheorem{lem}[thm]{Lemma}
\newtheorem{cor}[thm]{Corollary}
\newtheorem{claim}{Claim}
\newtheorem*{claim*}{Claim}

\newtheorem{mainThm}{Theorem}[]

\theoremstyle{definition}
\newtheorem{defn}[thm]{Definition}
\newtheorem{example}[thm]{Example}

\newtheorem*{examples*}{Examples}
\newtheorem*{example*}{Example}

\theoremstyle{remark}

\title[Quantum Homology of Real Lagrangians]{On the Quantum Homology
    of Real Lagrangians in Fano Toric Manifolds} 
\author{Luis Haug}
\thanks{The author was partially supported by the Swiss National
    Science Foundation.}
\address{Luis Haug, Department of Mathematics, ETH Z\"urich,
    R\"amistrasse 101, 8092 Z\"urich, Switzerland}
\email{haug@math.ethz.ch}

\date{\today}

\begin{document}
\begin{abstract}
    Let $R$ be the real Lagrangian in a toric symplectic manifold $X$,
    i.e. the fixed point set under complex conjugation. Assuming that
    $X$ is Fano of minimal Chern number at least 2, we show that the
    Lagrangian quantum homology $QH(R;\Lambda)$ of $R$ over the ring
    $\Lambda = \mathbb Z_2 [t, t^{-1}]$ is isomorphic to $H(R;\mathbb
    Z_2) \otimes \Lambda$ as a module. Moreover, we prove that $QH(R;
    \Lambda)$ is isomorphic \emph{as a ring} to the quantum homology
    $QH(X;\Lambda)$ of the ambient symplectic manifold.
\end{abstract}
\maketitle
\tableofcontents

\section{Introduction}
\label{sec:introduction}

\subsection{Setting}
\label{sec:intr-main-result}
Let $(X,\omega)$ be a symplectic toric manifold of dimension $2n$,
that is, a compact symplectic manifold endowed with an effective
Hamiltonian action of an $n$--dimensional torus $Q$. Denote by $\mu: X
\to \mathfrak q^*$ its moment map, where $\mathfrak q^*$ is the dual
of the Lie algebra of $Q$, and by $\Delta$ its moment polytope,
i.e., the image of $\mu$. $X$ carries a distinguished anti-symplectic
involution $\tau: X \to X$, characterized by the fact that it leaves
the moment map $\mu$ invariant. We refer to it as \emph{complex
    conjugation}. Its fixed point set
\begin{equation*}
    R \,:=\, \mathrm{Fix}\,\tau
\end{equation*}
is a Lagrangian submanifold, which we call the \emph{real Lagrangian}
of $X$. The prototypical example for a pair $(X,R)$ is $(\mathbb C
P^n, \mathbb R P^n)$.

For the purpose of this paper, we think of a symplectic toric manifold
as being obtained by equipping a Fano toric manifold $X$ with a
K\"ahler form $\omega$ whose cohomology class represents a positive
multiple of $c_1(X)$. In particular, we thus assume $\omega$ to be
monotone, which by definition means that the homomorphisms $I_\omega,
I_{c_1}: \pi_2(X) \to \mathbb R$ given by integration of the classes
$\omega$ and $c_1(X)$ satisfy
\begin{eqnarray*}
    I_{c_1} = \rho \,I_{\omega}
\end{eqnarray*}
for some $\rho \geq 0$. This ensures that also $R$ is monotone,
meaning in this case that the homomorphisms $I_{\omega}, I_{\mu}:
\pi_2(X,R) \to \mathbb R $ given by integration of $\omega$
respectively the Maslov index are related by
\begin{eqnarray*}
    I_{\omega} = \rho'\, I_{\mu}
\end{eqnarray*}
for some $\rho' \geq 0$. The minimal Chern number of $X$ and the
minimal Maslov number of $R$ are defined by
\begin{eqnarray*}
    C_X &:=& \min \big\{ \langle c_1(X), \lambda \rangle > 0 ~\big\vert~
    \lambda \in \pi_2(X) \big\},\\
    N_R &:=& \min \big\{ \langle \mu(R), \lambda \rangle > 0 ~\big \vert~
    \lambda \in \pi_2(X,R) \big\}.
\end{eqnarray*}
It is easy to see that they coincide in our situation, $C_X = N_R$. We
assume further that
\begin{equation*}
    C_X, N_R \geq 2.
\end{equation*}

\subsection{Main result}
\label{sec:main-result}

Under these assumptions, various versions of the \emph{Lagrangian
    quantum homology ring} $QH(R)$ introduced by Biran and Cornea
\cite{biran-cornea_ALagQH09} are defined. The purpose of this paper is
to study
\begin{equation*}
    QH(R; \Lambda_R),
\end{equation*}
the version obtained by using as coefficients the ring $\Lambda_R =
\mathbb Z_2[t^{-1},t]$ of $\mathbb Z_2$--Laurent polynomials, the latter
graded by setting the degree of the formal variable to be $\vert t
\vert = - N_R$. We will show that there is a close relation to the
quantum homology ring
\begin{equation*}
    QH(X; \Lambda_X)
\end{equation*}
of the ambient manifold defined over the ring $\Lambda_X = \mathbb
Z_2[q^{-1},q]$, where in this case the formal variable has degree
$\vert q \vert = - 2C_X$.

The main result of the paper is the following theorem.
\begin{mainThm}
    \label{Thm::Wideness_of_R}
    Let $(X,\omega)$ be a monotone symplectic toric manifold with
    minimal Chern number $C_X \geq 2$ and let $R$ be its real
    Lagrangian. Then the following hold.
    \begin{enumerate}
        \renewcommand{\theenumi}{\roman{enumi}}
    \item \label{Thm-part::wideness}There exists a canonical
        isomorphism of $\Lambda_R$--modules
        \begin{equation*}
            QH(R; \Lambda_R)  \cong  H(R; \mathbb Z_2) \otimes \Lambda_R.
        \end{equation*}
    \item \label{Thm-part::ring_str} Moreover, there exists a
        canonical degree-doubling ring isomorphism
        \begin{equation*}
            QH(R;\Lambda_R) \cong  QH(X;\Lambda_X).
        \end{equation*} 
    \end{enumerate}
\end{mainThm}

The first statement says that $R$ is $\Lambda_R$\emph{--wide} in the
terminology of \cite{biran-cornea_ALagQH09}, meaning the existence of
an isomorphism as in (i) (in general not required to be
canonical). Since $QH(R;\Lambda_R)$ is isomorphic to the Lagrangian
intersection Floer homology $HF(R;\Lambda_R)$, see
\cite{biran-cornea_ALagQH09}, the wideness immediately implies the
following non-displacement result: For any Hamiltonian diffeomorphism
$\phi$ of $(X,\omega)$ such that $R$ and $\phi(R)$ intersect
transversely, we have
\begin{equation*}
    \# \big(R \cap \phi(R)\big) \,\geq\, \sum_{j=1}^n b_j(R;\mathbb Z_2),
\end{equation*}
where the $b_j(R;\mathbb Z_2)$ denote the $\mathbb Z_2$ Betti numbers
of $R$. In terms of the moment polytope $\Delta$, the right-hand side
is the number of its vertices, as will become clear later.

The second part of the theorem allows to write down the ring structure
of $QH(R;\Lambda_R)$ using the known description of $QH(X;\Lambda_X)$
for Fano toric manifolds. The latter was first computed by Batyrev
\cite{Batyrev--QH_rings_of_toric_mfds__1993} and later by
Cieliebak--Salamon \cite{Cieliebak-Salamon--Wall_crossing__2006} using
different methods (both work with more general coefficients). We give
two examples.

\begin{example} Consider $X = \mathbb C P^n$ with its standard
    monotone symplectic form and its real part $R = \mathbb R P^n$. We
    have $C_X = N_R = n+1$. Denote by $h_{\mathbb C} = [\mathbb C P^{n-1}] \in
    H_{2n-2}(\mathbb C P^n;\mathbb Z_2)$ and $h_{\mathbb R} = [\mathbb
    R P^{n-1}] \in H_{n-1}(\mathbb R P^n;\mathbb Z_2)$ the classes of
    a complex respectively real projective hyperplane. These generate
    the classical intersection products: 
    $h_{\mathbb C}^{\cap i} = [\mathbb C P^{2n-2i}] \in H_{2n-2i}(\mathbb C
    P^n;\mathbb Z_2)$ and $h_{\mathbb R}^{\cap i} = [\mathbb R
    P^{n-i}] \in H_{n-i}(\mathbb R P^n;\mathbb Z_2)$. The quantum
    product on $QH(\mathbb C P^n;\Lambda_X) \cong H(\mathbb C
    P^n;\mathbb Z_2) \otimes \Lambda_X$ is given by
    \begin{equation*}
            h_{\mathbb C}^{\ast i} = 
        \begin{cases}
            h_{\mathbb C}^{\cap i},&i =1,\dots,n,\\
            [\mathbb C P^n] q, &  i = n+1,
        \end{cases}
    \end{equation*}
    see e.g. \cite{mcduff-salamon04}. It follows that the quantum ring
    structure on $QH(\mathbb R P^n;\Lambda_R) \cong H(\mathbb R
    P^n;\mathbb Z_2) \otimes \Lambda_R$ is given by
    \begin{equation*}
            h_{\mathbb R}^{\ast i} = 
        \begin{cases}
            h_{\mathbb R}^{\cap i},&i =1,\dots,n,\\
            [\mathbb R P^n] t, & i = n+1.
        \end{cases}
    \end{equation*}
\end{example}

\begin{example} Consider now $X = \mathbb C P^1 \times \mathbb C
    P^1\cong S^2 \times S^2$
    with the monotone symplectic form $\omega_{\mathbb C P^1} \times
    \omega_{\mathbb C P^1}$. The real part is $R = \mathbb R P^1
    \times \mathbb R P^1 \cong T^2$ and we have $C_X = N_R =
    2$. Setting $A := [\mathbb C P^1 \times \mathrm{pt}]$,
    $B:=[\mathrm{pt} \times \mathbb C P^1] \in H_2(X;\mathbb Z_2)$,
    the quantum product on $QH(X;\Lambda_X) = H(X;\mathbb Z_2)
    \otimes \Lambda_X$ satisfies
    \begin{equation*}
        A \ast B = B \ast A = [\mathrm{pt}], \quad A \ast A = B \ast B =
        [X] q.
    \end{equation*}
    (In fact, this determines all other products by the associativity of
    $\ast$.)  Setting $a := [\mathbb R P^1 \times \mathrm{pt} ], b:=
    [\mathrm{pt} \times \mathbb R P^1] \in H_1(R;\mathbb Z_2)$, we
    conclude that the quantum product on $QH(R;\Lambda_R) \cong
    H(R; \mathbb Z_2) \otimes \Lambda_R$ is determined by
    \begin{equation*}
        a \ast b = b \ast a = [\mathrm{pt}], \quad a \ast a = b \ast b
        =  [R]t.
    \end{equation*}
\end{example}

We now give a heuristic explanation for why the wideness part of
Theorem \ref{Thm::Wideness_of_R} should hold. The quantum homology
$QH(R;\Lambda_R)$ is the homology of a deformed Morse
complex, whose differential $d$ decomposes as
\begin{equation*}
    d = d_0 + d_1.
\end{equation*}
The ``classical'' part $d_0$ counts Morse trajectories, while the
``quantum'' part $d_1$ counts pearly trajectories, that is, certain
configurations of Morse trajectories with (pseudo)--holomorphic discs
interspersed between them. Whenever we work with a
$\tau$--anti-invariant almost complex structure, such as the standard
integrable $J_0$, there exists an involution on the moduli spaces of
discs, given by reflecting discs along $R$, and consequently also on
the moduli spaces of pearly trajectories. The idea is that all the
quantum contributions entering $d_1$ should cancel in the $\mathbb
Z_2$--count because of this involution, so that $d_1$ vanishes and
hence $d$ is equal to the Morse differential. However, this reasoning
does not take into account that the involution might have fixed points
(it actually does!), which a priori destroys the argument. Showing
that it still works requires to take a closer look at the involution,
see Section \ref{sec:invol-moduli-spac}.

To approach the relation between $QH(R;\Lambda_R)$ and
$QH(X;\Lambda_X)$, the first thing to note is that on the classical
level there exists a degree-doubling isomorphism
\begin{equation*}
    H(R;\mathbb Z_2) \to H(X;\mathbb Z_2)
\end{equation*}
which respects the ring structures given by the intersection
products. The inverse map is easy to describe: It is induced by
mapping all $\tau$--invariant cycles, which generate $H(X;\mathbb
Z_2)$, to their $\tau$--fixed parts. We will use an argument of
Duistermaat, ultimately based on classical Smith theory, to show that
this leads to an isomorphism. While it seems unlikely that this
isomorphism has not been known before, the author is not aware of a
position in the literature where it appears.

Part (\ref{Thm-part::ring_str}) of Theorem \ref{Thm::Wideness_of_R}
states that this relation stays true on the quantum level. The
heuristic idea is similar to the one for the proof of the wideness
statement. Since $QH(X;\Lambda_X)$ is by definition $H(X;\mathbb
Z_2) \otimes \Lambda_X$ as a $\Lambda_X$--module, the wideness of $R$
allows to conclude immediately that $QH(R;\Lambda_R)$ and
$QH(X;\Lambda_X)$ are isomorphic as modules. In the case of the
ambient manifold $X$, the quantum product is a deformation of the
classical intersection product using 3--point genus zero
Gromov--Witten invariants. In the Morse picture, this means that the
non-classical part counts Y--type configurations of Morse trajectories
connected by a holomorphic sphere in the center. In the case of the
real part $R$, the quantum product counts Y--type configurations of
pearly trajectories connected by a holomorphic disc. Similar to the
case of the differential, there exists an involution on the moduli
space of Y--pearls, and the same type of cancellation argument says
that all contributions with a non-central non-constant disc should
occur in pairs and hence cancel when $\mathbb Z_2$--counting
them. Moreover, if we work with $\tau$--invariant Morse functions $f:
X \to \mathbb R $ having all their critical points on $R$, there
should be a 1--1 correspondence between the remaining Y--pearls on the
Lagrangian side and those on the ambient side.

Again, this argument fails a priori because the involution on the
moduli space of discs has fixed points.

\subsection{Previous results} 
\label{sec:context}
The first results in the direction of Theorem \ref{Thm::Wideness_of_R}
were obtained by Oh \cite{Oh--HF_of_Lagr_Int_II--1993} who proved that
$HF(\mathbb R P^n;\mathbb Z_2)$ is isomorphic as a vector space to
$H(\mathbb R P^n; \mathbb Z_2)$. Later Biran and Cornea computed the
ring structure of $QH(\mathbb R P^n; \Lambda)$, in
fact as a special case of a series of more general results for
Lagrangians $L \subset \mathbb C P^n$ satisfying the condition
$2H_1(L;\mathbb Z) = 0$ (see \cite{biran-cornea_rigidity09} or Section
6.2.2 in \cite{biran-cornea_QuantumStrForLag07}).

As for general toric manifolds, there is another family of natural
Langrangian submanifolds, namely the smooth fibres of the moment
map. The study of the Floer homology of the monotone ones was started
by Cho--Oh \cite{Cho-Oh_FloerCohomology2006} and later extended in a
series of papers by Cho. There are also recent results due to
Alston--Amorim \cite{Alstom-Amorim--HF_torus_real--2010} and
Abreu--Macarini \cite{Abreu-Macarini--Rmks_Lagr_Intersections--2011}
on the Floer homology $HF(R,L)$ of the real part $R$ and a smooth
moment fibre $L$.

There have also been efforts to calculate the Floer homology (as a
module) of Lagrangians which are fixed point sets of anti-symplectic
involutions in more general contexts than that of toric manifolds, see
e.g. the works by Frauenfelder
\cite{Frauenfelder--AG_conjecture--2004} and Fukaya--Oh--Ohta--Ono
\cite{FOOO--Anti-symplectic_involution_and_HF--2009}.

\subsection{Organisation of the paper}
\label{sec:outline-paper}
We provide the necessary background on symplectic toric manifolds in
Section \ref{sec:sympl-toric-manif}. In Section
\ref{sec:topology-real-part} we examine the topology of $R$ and its
relation to that of $X$ using $\tau$--invariant Morse functions. More
precisely, we prove $H(R;\mathbb Z_2) \, \cong \, H(X;\mathbb Z_2)$ as
rings. Section \ref{sec:real-discs-real} provides a sort of
classification of the holomorphic discs $u: (D^2,\partial D^2) \to (X,R)$
with boundary on $R$. We show that after after removing a boundary
point, they admit a representation in toric homogeneous
coordinates. This is used in Section \ref{sec:transversality} to prove
that the standard complex structure $J_0$ is regular for
them. Basically the same argument will also show the regularity of
$J_0$ for all holomorphic spheres $v: \mathbb C P^1 \to X$. In Section
\ref{sec:invol-moduli-spac} we take a closer look at the involution on
the moduli space of real discs. We describe its fixed point set and
show that all fixed points we eventually have to count in the
definition of the quantum invariants occur in pairs. Section
\ref{sec:transv-eval-maps} is technical and serves to show that
generic choices of $\tau$--invariant data lead to well-behaved moduli
spaces. In Section \ref{sec:proofs-main-thms} we finally assemble all
the ingredients and prove the main result.

\subsection* {Acknowledgement}
\label{sec:acknowledgement}
The author wishes to thank his advisor, Professor Paul Biran, for
suggesting to start this investigation and for plenty of helpful
discussions on the subject. He also wishes to thank the anonymous
referee for his comments.

\section{Symplectic Toric Manifolds}
\label{sec:sympl-toric-manif}
There are essentially two ways of constructing a symplectic toric
manifold, both starting from a suitable polytope $\Delta$ in $(\mathbb
R^n)^*$. The first is by symplectic reduction of some $\mathbb C^N$
with respect to a torus action, in a way encoded by $\Delta$. In the
second, one starts by building a complex manifold $X_\Sigma$ together
with a torus action from the normal fan $\Sigma$ of $\Delta$, and then
embeds $X_\Sigma$ into some big projective space to define the
symplectic form. In both cases one obtains a symplectic manifold
together with a Hamiltonian torus action whose moment polytope is
$\Delta$.

The two constructions are equivalent in the sense that there exist
equivariant symplectomorphisms between the resulting symplectic toric
manifolds, which follows from a theorem of Delzant
\cite{Delzant--Hamiltoniens_periodiques__1988}. We will describe the
second variant as it is more suited to our needs.

The entire section is heavily inspired by the texts from which
the author first learned about symplectic toric manifolds, mainly the
book by Audin \cite{Audin_TorusActions2004} and the lecture notes by
Cannas da Silva contained in
\cite{Audin-Cannas_da_Silva-Lerman__2003}, but also to some extent the
treatment in Cho and Oh's paper \cite{Cho-Oh_FloerCohomology2006}.

\subsection{Cones, fans, polytopes}
\label{sec:cones-fans}
We start by giving definitions of the relevant combinatorial objects.

\begin{defn}
    A \emph{cone} in $\mathbb R^n$ of dimension $r$ is a subset
    $\sigma \subset \mathbb R^n$ that can be written in the form
    \begin{equation*}
        \sigma = \{ a_1 v_1 + \dots + a_r v_r ~|~
        a_1,\dots,a_r \geq 0\}
    \end{equation*}
    with linearly independent vectors $v_1,\dots, v_r \in \mathbb
    R^n$, the \emph{generators} of $\sigma$. A cone is called
    \emph{rational and smooth} if it admits a set of generators in
    $\mathbb Z^n$ which can be extended to a $\mathbb Z$--basis of
    $\mathbb Z^n$. Any other cone $\sigma' \subset \mathbb R^n$
    admitting a set of generators $v_1',\dots,v_s'$ such that $\{v_1',
    \dots, v_s' \} \subset \{v_1,\dots,v_r\}$ is called a \emph{face}
    of $\sigma$.
\end{defn}

\begin{defn}
    A \emph{fan} in $\mathbb R^n$ is a finite set $\Sigma =
    \{\sigma_1,\dots, \sigma_\ell\}$ of cones in $\mathbb R^n$ with the
    following properties: \renewcommand{\theenumi}{\roman{enumi}}
    \begin{enumerate}
    \item For every cone $\sigma \in \Sigma$ and every face $\sigma'$
        of $\sigma$, we have $\sigma' \in \Sigma$.
    \item For all cones $\sigma, \sigma' \in \Sigma$, the cone $\sigma
        \cap \sigma'$ is a face of both $\sigma$ and $\sigma'$.
    \end{enumerate}
    $\Sigma$ is \emph{complete} if in addition the condition
    \begin{enumerate}
    \setcounter{enumi}{2}
    \item $\sigma_1 \cup \dots \cup \sigma_\ell = \mathbb R^n$
    \end{enumerate}
    holds. For a fan $\Sigma$, we denote by $\Sigma^{(r)}$ the set of
    its $r$-dimensional cones.
\end{defn}

\begin{defn}
\label{defn::polytope}
    A \emph{polytope} $\Delta$ in $(\mathbb R^n)^*$ is a bounded
    subset that can be written as the intersection of a finite number
    of affine half-spaces, i.e.
    \begin{equation*}
        \Delta \,=\, \big\{ \varphi \in (\mathbb R^n)^* ~\big\vert~ \langle
        \varphi, v_i \rangle \geq - a_i  ~\forall ~i \in [1,N] \big\}
    \end{equation*}
    with $v_1,\dots,v_N \in \mathbb R^n$ and $a_1,\dots, a_N \in
    \mathbb R$. Given a fixed collection of such defining normal
    vectors $v_i$, for $I \subset [1,N]$ the subset
    \begin{equation*}
       F_I \,:=\, \big\{
        \varphi \in \Delta ~|~ \langle \varphi, v_i \rangle = - a_i
        ~\forall~ i \in I \big \}
    \end{equation*}
    is called the \emph{$I$-th face} of $\Delta$ if it is non-empty.
\end{defn}

A polytope $\Delta$ in $(\mathbb R^n)^*$ gives rise to a fan $\Sigma
(\Delta)$ in $\mathbb R^n$ as follows:
\begin{defn}
    Given a face $F_I$ of $\Delta$, let $\sigma_I$ be the cone
    generated by the normal vectors $v_i$, $i \in I$. Then
    \begin{equation*}
        \Sigma(\Delta) \, := \, \big\{ \sigma_I ~\big\vert~ \text{ $F_I$ is a face
            of $\Delta$} \big\}
    \end{equation*}
    is a fan called the \emph{normal fan} of $\Delta$.
\end{defn}

\subsection{Constructing a toric manifold from a fan}
\label{sec:constr-toric-manif}
Let $\Sigma$ be a complete fan of smooth rational cones of dimension
$n$. In this subsection we will construct the \emph{toric manifold}
associated to $\Sigma$, which by definition is a compact complex
manifold
\begin{equation*}
    X \equiv X_\Sigma
\end{equation*}
together with an effective action of an $n$--dimensional complex
torus $Q_{\mathbb C} \cong T_{\mathbb C}^n.$\footnote{We write $T_{\mathbb C}^n$ for $(\mathbb C^*)^n$
    whenever it plays the role of a group acting on a space.}

The starting point for the construction of $X$ is the standard action
of $T_{\mathbb C}^N$ on $\mathbb C^N$ given by
\begin{equation*}
   (t_1,\dots,t_N) \cdot (z_1,\dots,z_N) = (t_1z_1,\dots,t_Nz_N).
\end{equation*}
The idea is to use the combinatorics of $\Sigma$ to find a certain
subtorus $K_{\mathbb C} \subset T_{\mathbb C}^N$ and a subset
$\mathcal U \subset \mathbb C^N$ such that the restriction of the
action of $K_{\mathbb C}$ to $\mathcal U$ is free. 

As for the definition of $K_{\mathbb C}$, denote by $\{v_1,\dots,
v_N\} \subset \mathbb Z^n$ the set of primitive integral generators of
the $1$--dimensional cones in $\Sigma$ and consider the homomorphism
\begin{equation*}
    \pi: \mathbb Z^N \to \mathbb Z^n, \qquad e_j \mapsto v_j \quad \text{for~}
    j = 1,\dots,N,
\end{equation*}
where $e_j = (0,\dots,1,\dots,0) \in \mathbb Z^N$ is the $j$--th
standard basis vector. Denote by $\mathbb K$ the kernel of
$\pi$ and define $K_{\mathbb C}$ to be the subgroup of
$T_{\mathbb C}^N$ generated by all 1--parameter subgroups $T_{\mathbb
    C}^1 \to T_{\mathbb C}^N$ of the form
\begin{equation*}
    \tau \,\mapsto\, \tau^\lambda \,:=\, (\tau^{\lambda_1},\dots,\tau^{\lambda_N}),
\end{equation*}
with $\lambda = (\lambda_1,\dots,\lambda_N) \in \mathbb
K$. 

Next we construct $\mathcal U \subset \mathbb C^N$. Given a subset $I
= \{i_1,\dots,i_r\} \subseteq [1,N]$, consider the linear
subspace $Z_I \,:=\, \{ z \in \mathbb C^N ~|~ i \in I \,\Longrightarrow \,
    z_i = 0 \}$
and denote by $\sigma_I$ the cone generated by the $v_i$, $i \in
I$. Then set
\begin{equation*}
    \mathcal U \, := \, \mathbb C^N \smallsetminus \bigcup_{\sigma_I
        \notin \Sigma} Z_I.
\end{equation*}
To describe this in words: A point $z = (z_1,\dots,z_N) \in \mathbb
C^N$ belongs to $\mathcal U$ iff the cone $\sigma_{I_z}$ belongs to
$\Sigma$, where $I_z := \big\{i \in [1,N]~\big\vert~ z_i=0\big\}$.

\begin{lem}
    \label{lem::freeness_of_KC-action}
    The restriction of the standard action of the subgroup $K_{\mathbb
        C} \subset T_{\mathbb C}^N$ to the set $\mathcal U \subset
    \mathbb C^N$ is free.
\end{lem}

It follows that the quotient
\begin{equation*}
    X \,\equiv X_\Sigma := \, \mathcal U / K_{\mathbb C}
\end{equation*}
is a complex manifold of dimension $n$, with the complex structure
$J_0$ inherited from the one on $\mathbb C^N$. We denote by $\pi:
\mathcal U \to X$ the canonical projection and commonly refer to points
of $X$ using ``toric homogeneous coordinates'', i.e. we write
\begin{equation*} 
    [z_1:\dots:z_N]_{X} \,:= \, \pi \big(z_1,\dots,z_N\big),
\end{equation*}
in reminiscence of the homogeneous coordinates on projective space.
The completeness of $\Sigma$ implies that $X$ is compact (see
Proposition VII.1.4 in \cite{Audin_TorusActions2004}). Moreover, the
$T_{\mathbb C}^N$--action on $\mathcal U$ induces an action of the
quotient
\begin{equation*}
    Q_{\mathbb C} \, := \, T_{\mathbb C}^N / K_{\mathbb C}
\end{equation*}
on $X$, which is an $n$--dimensional complex torus.

\subsection{The toric divisors, the first Chern class and the homology
    ring}
\label{sec:toric-divisors}
To any cone $\sigma \in \Sigma$ of dimension $r$, we associate a
complex submanifold of $X$ by setting 
\begin{equation*}
    C_\sigma \, := \, \pi (Z_{I_\sigma}),
\end{equation*}
where $I_\sigma := \big\{i \in [1,N] ~\big\vert~ \text{$v_i$ generates
    an edge of $\sigma$}\big\}$. $C_\sigma$ has complex codimension
$r$ and is invariant under the $Q_{\mathbb C}$--action because
$Z_{I_\sigma}$ is invariant under the standard $T_{\mathbb
    C}^N$--action. Using toric homogeneous coordinates, we can write
\begin{equation*}
    C_\sigma \, = \, \big\{ [z_1: \dots : z_N]_{X} ~\big\vert~
    i \in I_{\sigma} \Longrightarrow z_i = 0 \big\}.
\end{equation*}

Of particular importance are the hypersurfaces corresponding to the
1-dimen\-sional cones in $\Sigma$. We denote them by
\begin{equation*}
    D_1,\dots,D_N \subset X
\end{equation*}
and refer to them as the \emph{toric divisors}. Their homology classes
generate $H_*(X;\mathbb Z_2)$ as a ring with respect to the
intersection product. More precisely, we have
\begin{equation*}
    H_*(X;\mathbb Z_2) \cong \mathbb Z_2 \big[ [D_1], \dots, [D_N] \big] \Big/
    (\mathcal I + \mathcal J),
\end{equation*}
where $\mathbb Z_2 \big[ [D_1], \dots, [D_N] \big]$ is the polynomial
ring in the variables $[D_i]$, and the ideals by which one quotients
out are
\begin{equation*}
    \mathcal I = \left\langle \sum_{i = 1}^N \langle \varepsilon, v_i\rangle [D_i] ~\Big|~ \varepsilon \in
        (\mathbb Z^n)^* \right\rangle,
    ~ \mathcal J = \left \langle \prod_{i \in I} [D_i] ~\Big|~ I \subset
        [1,N] \text{~with~} \sigma_{I} \notin \Sigma \right \rangle.
\end{equation*}
Moreover, the union of the $D_i$ represents the first Chern class of
$X$, in the sense that
\begin{equation*}
    c_1(X) \,=\, PD \big ([D_1] + \dots + [D_N] \big).
\end{equation*}
For all of this consult \cite{Audin_TorusActions2004}.

\subsection{Affine open subsets}
\label{sec:toric-divisors}
We next introduce a cover of $X$ by open subsets biholomorphic to
$\mathbb C^n$ corresponding to the $n$--dimensional cones in
$\Sigma$. Given any such cone $\sigma \in \Sigma^{(n)}$, let $I_\sigma
= \big\{i \in [1,N] ~\big\vert~ \text{$v_i$ generates an edge of
    $\sigma$}\big\}$ be as before and set
\begin{equation*}
    V_\sigma \,:= \, \big\{ [z_1 : \dots : z_N]_{X} \in
    X  ~\big\vert~ z_i = 0 ~\Longrightarrow~
    i \in I_\sigma \big\}.
\end{equation*}
That is, $V_\sigma$ is obtained from $X$ by cutting out those
toric divisors $D_i$ for which $v_i$ is not part of the set of
generators of $\sigma$. 

It is easy to see that the $V_\sigma$ with $\sigma \in \Sigma^{(n)}$
provide a cover of $X$. Indeed, given any point $z \in \mathcal U$,
let $I_z \,:= \, \big\{i \in [1,N] ~\big\vert~ z_i = 0 \big\}$ and
note that we have $\sigma_{I_z} \in \Sigma$ by the definition of
$\mathcal U$. Hence $\sigma_{I_z}$ is a face of an $n$--dimensional
cone $\sigma \in \Sigma$ by the completeness of $\Sigma$, and clearly
we have $z \in V_\sigma$.

\begin{lem}
\label{lem::affine_open_sets}
 For every $\sigma \in \Sigma^{(n)}$, the set $V_\sigma
    \subset X$ is biholomorphic to $\mathbb C^n$.
\end{lem}
\begin{proof}
    We assume without loss of generality that the generators of the
    edges of $\sigma$ are the vectors $v_{k + 1},\dots,v_N \in \mathbb
    Z^n $, where $k := N-n$. These form a basis of $\mathbb Z^n$, and
    we denote by $\nu_{k+1},\dots,\nu_N \in (\mathbb Z^n)^*$ the dual
    basis.
    
    Next we define a basis $\lambda^{(1)},\dots,\lambda^{(k)} \in
    \mathbb Z^N$ of $\mathbb K = \mathrm{ker}\,\pi$ by setting
    \begin{equation*}
        \lambda^{(i)}_j \,:=\,
        \begin{cases}
            \delta_{ij}, & j \in [1,k]\\
            -\langle \nu_j,v_i \rangle, & j \in [k+1,N].
        \end{cases}
    \end{equation*}
    Note that we have $\pi (\lambda^{(i)}) = v_i - \sum_{j=k+1}^N
    \langle \nu_j,v_i \rangle v_j$ for all $i \in [1,k]$, and
    consequently $\langle \nu_j, \pi (\lambda^{(i)}) \rangle = \langle
    \nu_j,v_i \rangle - \langle \nu_j,v_i\rangle = 0$ for $j \in
    [k+1,N]$, which shows that all $\lambda^{(i)}$ belong to $\mathbb
    K$. Since they are clearly linearly independent, they constitute
    indeed a basis of $\mathbb K$.

    Let now $z = [z_1:\dots:z_N]_{X} \in V_\sigma$, which means
    that $z_1,\dots,z_k \neq 0$. We claim that $z$ is represented by
    an element of $\mathcal U$ of the form $(1,\dots,1,
    \ast_{k+1},\dots,\ast_N)$. Namely, acting on such a vector with
    $z_1^{\lambda^{(1)}} \dotsm z_k^{\lambda^{(k)}} \in K_{\mathbb C}$
    yields
    \begin{equation*}
            z_1^{\lambda^{(1)}} \dotsm z_k^{\lambda^{(k)}} \cdot
            (1,\dots,1,\ast_{k+1},\dots,\ast_N) \,=\,
            (z_1,\dots,z_k,\ast_{k+1}',\dots,\ast_N')
    \end{equation*}
    with $\textstyle \ast_j' \,=\, \ast_j \cdot \prod_{i=1}^k
    z_i^{\lambda_j^{(i)}}$ for $j \in [k+1,N]$. Setting $\ast_j
    \,:=\, z_j \cdot \prod_{i=1}^k z_i^{-\lambda_j^{(i)}}$ thus
    produces the required vector. Note that we can write
    \begin{eqnarray*}
        \textstyle
        \ast_j \,=\, z_j \cdot \prod_{i = 1}^k z_i^{\langle \nu_j,v_i
            \rangle} \,=\,  \prod_{i=1}^N z_i^{\langle \nu_j,v_i\rangle},
    \end{eqnarray*}
    because $\langle \nu_j, v_i \rangle = \delta_{ij}$ for $i \in
    [k+1,N]$. It follows that the map given by
    \begin{equation*}
        \textstyle
        [z_1: \dots: z_N]_{X} \,\mapsto\, \left( \prod_{i=1}^N
            z_i^{\langle \nu_{k+1},v_i \rangle},\dots, \prod_{i=1}^N
            z_i^{\langle \nu_{N},v_i \rangle} \right)
    \end{equation*}
    provides a biholomorphism $V_\sigma \,\cong\, \mathbb C^n$.
\end{proof}

\subsection{The symplectic form determined by a polytope}
\label{sec:sympl-form-determ}
Let $\Delta$ be a polytope in $(\mathbb R^n)^*$ and let
$\Sigma(\Delta)$ be its normal fan. In this subsection we describe how
$\Delta$ determines an embedding of the toric manifold $ X \equiv
X_{\Sigma(\Delta)}$ into some projective space. Pulling back the
Fubini-Study form $\omega_{FS}$ will then equip $X$ with a symplectic
form $\omega$ such that the action of the real torus $Q \subset
Q_{\mathbb C}$ is Hamiltonian. For this to work we need to assume that
$\Delta$ has the following additional properties:
\begin{enumerate}
\item $\Delta$ is \emph{lattice}, in the sense that all vertices
    belong to $(\mathbb Z^n)^*$.

\item $\Delta$ is \emph{Delzant}, meaning that its defining set of
    normal vectors $v_i$, $i \in [1,N]$, can be chosen such that
    for every (non-empty) face $F_I$, the corresponding collection of
    $v_i$, $i \in I$, can be extended to a $\mathbb Z$--basis of
    $(\mathbb Z^n)^*$.
\end{enumerate}
$\Sigma(\Delta)$ is then a complete fan of smooth rational cones, and
hence $X$ is a compact toric manifold.

The construction of the embedding goes as follows. Denote by
$m_1,\dots, m_L \in \Delta \cap (\mathbb Z^n)^*$ the lattice points of
$\Delta$, and let $F_i \subset \Delta$ be the $i$--th facet of
$\Delta$, i.e. the codimension--1 face with normal vector $v_i$. Set
$LD(m_i, F_j) \,:=\, \langle m_i, v_j \rangle + a_j \in \mathbb Z$,
which is the lattice distance of $m_i$ from the facet $F_j$. That is,
$m_i$ lies on the $LD(m_i,F_j)$--th hyperplane parallel to $F_j$,
counting from $F_j$. Now define a map $\mathbb C^N \to \mathbb C^L$ by
\begin{equation*}
    (z_1,\dots,z_N) \, \mapsto\, \left(\prod_{j=1}^N
        z_j^{LD(m_1,F_j)},\dots,\prod_{j=1}^N
        z_j^{LD(m_L,F_j)}  \right).
\end{equation*}
Its restriction to $\mathcal U \subset \mathbb C^N$ takes values in
$\mathbb C^L \smallsetminus \{0\}$ and descends to a map
\begin{equation*}
    \Psi: X \, \to \, \mathbb C P^{L-1}
\end{equation*}
which is actually an embedding (see
e.g. \cite{Oda--Convex_bodies--1987}). Moreover, it is equivariant
with respect to the action of the real torus
$Q := \big\{(t_1,\dots,t_N) \in Q_{\mathbb C} ~\big\vert~ \vert t_i \vert
    = 1 ~\forall \,i \big\}$
on $X$ and the $\mathbb T^n$--action on $\mathbb C
P^{L-1}$ given by
\begin{equation*}
    (t_1,\dots,t_n) \cdot \big[z_1:\dots:z_L\big] = \big[t^{m_1}
    z_1: \dots: t^{m_L} z_L\big],
\end{equation*}
with $t ^{m_i} := t_1^{m_{i,1}}\cdot \dots \cdot t_n^{m_{i,n}}$. The
latter is Hamiltonian and has moment map
\begin{equation*}
    \hat\mu: \mathbb C P^{L-1} \to \mathbb R^n, \quad [z_1:\dots:z_L]
    \mapsto \frac{\sum_{i=1}^L m_i \vert z_i \vert^2}{\sum_{i=1}^L
        \vert z_i \vert^2},
\end{equation*}
with respect to (a multiple of) $\omega_{FS}$. The image of $\hat \mu$
is the convex hull of the $m_i$, i.e. the polytope $\Delta$, and this
stays true for the restriction of $\hat \mu$ to $\Psi(X)$. Consult the
last section of Cannas da Silva's notes
\cite{Audin-Cannas_da_Silva-Lerman__2003} for more details on that.

It follows that $\omega := \Psi^* \omega_{FS}$ makes $X$ a symplectic
manifold such that the $Q$--action is Hamiltonian with a moment map
$\mu: X \to \mathfrak q^*$ whose image is the polytope $\Delta
\subset \mathfrak q^* \cong (\mathbb R^n)^*$.

\subsection{Complex conjugation}
\label{sec:complex-conjugation}
Complex conjugation on $\mathcal U \subset \mathbb C^N$ descends to an
involution
\begin{equation*}
    \tau: X \to X
\end{equation*}
which we also refer to as complex conjugation. It is anti-symplectic
with respect to $\omega = \Psi^*\omega_{FS}$ and related to the torus
action by
\begin{equation*}
    \tau \circ q = q^{-1} \circ \tau
\end{equation*}
for $q \in Q$. In particular, $\tau$ maps $Q$--orbits to $Q$--orbits
but reverses the time order on them. Observe that $\tau$ leaves the
moment map invariant, i.e. $\mu \circ \tau = \mu$. This property
characterizes $\tau$ among all anti-symplectic involutions on $X$
(up to conjugation by $Q$--equivariant symplectomorphisms).

\section{The Topology of the Real Lagrangian}
\label{sec:topology-real-part}

Let now $X$ be a symplectic toric manifold of real dimension $2n$ and
$R$ its real part. The purpose of this section is to explain the
relation between the homology rings $H_*(X;\mathbb Z_2)$ and
$H_*(R;\mathbb Z_2)$.

\subsection{$\tau$--invariant Morse functions}
\label{sec:tau-invariant-morse}
Denote by $Q$ the $n$--dimensional torus acting on $X$ in a
Hamiltonian way, and by $\mu: X \to \mathfrak q^*$ the corresponding
moment map, where $\mathfrak q^*$ is the dual of the Lie algebra
$\mathfrak q$ of $Q$. We will work with the following family of
$\tau$--invariant functions, defined by
\begin{equation*}
    f_\xi: X \to \mathbb R, ~ x \mapsto \langle \mu(x), \xi \rangle,
\end{equation*}
for $\xi \in \mathfrak q$. If $\xi$ is chosen generically, then
$f_\xi$ is Morse, and in this case the critical points of $f_\xi$ are
precisely the points where $D \mu$ vanishes, i.e., the preimages of
the vertices of the moment polytope $\Delta$. The Morse index $\vert p
\vert$ at such a point is given by twice the number of edges of
$\Delta$ meeting at the vertex $\mu(p)$ for which the evaluation of an
inward-pointing generating vector on $\xi$ is negative. In particular,
all critical points have even Morse index, so that $f_\xi$ is perfect,
i.e. the Morse differential vanishes.

The restriction $f_\xi\vert_R$ is a Morse function on $R$, and the
$\tau$-invariance implies easily that all critical points of
$f_\xi\vert_R$ are also critical points of $f_\xi$,
i.e. $\mathrm{Crit} f_\xi\vert_R \subseteq \mathrm{Crit} f_\xi$. In
fact we even have:

\begin{lem}
    \label{lemma:Critical_points_on_X_and_R}
    The sets $\mathrm{Crit} f_\xi$ and $\mathrm{Crit} f_\xi \vert_R$
    coincide, and the index of any critical point $p$ of $f_\xi$ is
    twice the index of $p$ as a critical point of $f_\xi\vert_R$.
\end{lem}
\begin{proof}
    For the first statement, it remains to show that $\mathrm{Crit}
    f_\xi \subseteq \mathrm{Crit} f_\xi \vert_R$, which clearly
    follows if we can show $\mathrm{Crit} f_\xi \subset R$. Note that
    the $\mu$--preimages of all faces of $\Delta$, in particular of
    the vertices, are connected and invariant under $\tau$, because
    their preimages under $\pi: \mathcal U \to X$ are the sets $Z_I
    \subset \mathbb C^N$ (see Section \ref{sec:constr-toric-manif}),
    which are connected and invariant under complex conjugation. In
    particular, the preimage of every vertex is a single point in $R$,
    and hence $\mathrm{Crit}\,f_\xi$ is contained in $R$.

    The statement about the indices is an immediate consequence of the
    fact that at every critical point $x$ of $f_\xi$, there exist
    local coordinates $q_j,p_j$ centered at $x$, such that $f_\xi$ can
    locally be expressed as
    \begin{equation*}
        f_\xi = f_\xi(x) + \frac{1}{2} \sum_{j=1}^n \langle \eta_j,\xi\rangle \big(q_j^2 + p_j^2\big)
    \end{equation*}
    with certain $\eta_j \in \mathfrak q^*$, and such that $R$ is
    given by the vanishing of $p$ locally. This is part of Proposition
    2.2 in \cite{duistermaat--convexity_and_tightness__1983}. In these
    coordinates the Hessians are
    \begin{eqnarray*}
        \mathrm{Hess}\,f_\xi\vert_R (x) &=&
        \mathrm{diag}(\langle \eta_1, \xi\rangle,\dots,\langle \eta_n,
        \xi\rangle),\\
        \mathrm{Hess}\,f_\xi (x) &=&
        \mathrm{diag}(\langle \eta_1, \xi\rangle,\dots,\langle \eta_n,
        \xi\rangle,\langle \eta_1, \xi\rangle,\dots,\langle \eta_n,
        \xi\rangle),
    \end{eqnarray*}
    from which the claimed relation follows.
\end{proof}

Recall that the definition of Morse homology for a pair $(f,g)$ of a
Morse function and a Riemannian metric requires Morse--Smale
transversality, meaning that all (un)stable manifolds
$W_x^u(f,g)$, $W_y^s(f,g)$ intersect transversely. It is a standard
fact that for a fixed Morse function, a generic choice of Riemannian
metric yields a Morse--Smale pair, see e.g. Proposition 5.8 in
\cite{Hutchings--Notes_on_Morse_Homology_2002}. In our situation, we
are constrained by the fact that we want not only the Morse function,
but also the Riemannian metric to be $\tau$--invariant, such as to
fully exploit the symmetry. However, the standard result stays true:

\begin{prop}
    \label{prop:tau-invt-Morse-Smale}
    Let $f: X \to \mathbb R$ be a $\tau$--invariant Morse
    function. There exists a second category subset of the space
    $\mathcal G^\tau$ of $\tau$-invariant Riemannian metrics on $X$
    such that for every $g$ in this subset the pair $(f, g)$ is
    Morse-Smale.
\end{prop}

The proof is almost identical to the proof of the statement without
the invariance condition. The only point where a bit of caution is
necessary is when local perturbations of the metric along negative
gradient flow lines have to be extended to global
$\tau$-\emph{invariant} perturbations. But this poses no problem,
basically because if $\gamma$ is such a flow line, we either have
$\mathrm{im}\,\gamma \subset R$ or $\mathrm{im}\,\gamma \cap
\mathrm{im}\,\tau\circ\gamma = \emptyset$, as a consequence of the fact
that both $f$ and $g$ are $\tau$--invariant. (See also Proposition
\ref{prop:regularity-mod-spaces} in this paper and its proof.)

\begin{lem} Let $f:X \to \mathbb R$ be a $\tau$--invariant Morse
    function and $g$ a $\tau$--invariant metric on $X$ such that
    $(f,g)$ is Morse-Smale. Then the restricted pair
    $(f\vert_R,g\vert_R)$ is also Morse-Smale.
\end{lem}

\begin{proof} The invariance implies that $-\nabla f$ is tangent to
    $R$. Thus the flow of $-\nabla f \vert_R$ is simply the
    restriction of the flow of $-\nabla f$ to $R$, and the
    invariant manifolds are related by $W_{x}^*(f\vert_R) = W_x^*(f)
    \cap R$, for $x \in \mathrm{Crit}\,f \cap R$. Given some $p \in R$,
    any subspace $V \subset T_pX$ with $\tau V = V$ decomposes as $V =
    V^+ \oplus V^-$, with $V^{\pm} := (\mathrm{id} \pm D\tau)V$
    denoting its $\tau$--invariant respectively $\tau$--anti-invariant
    subspace, which are both of the same dimension $\frac{1}{2}\,
    \mathrm{dim}\,V$. Now $T_p W_x^*(f\vert_R) = T_p(W_x^*(f) \cap R)$
    is clearly contained in $(T_pW_x^*(f))^+ = T_p W_x^*(f) \cap
    T_pR$, but since both have dimension $\frac{1}{2} \mathrm{dim}\,
    W_x^*(f)$, we actually have $T_pW_x^*(f\vert_R) = (T_p
    W_x^*(f))^+$. (In other words, all $W_x^*(f)$ intersect $R$
    cleanly).

    Now $(f,g)$ being Morse-Smale means that $T_p W_x^u(f) + T_p
    W_y^s(f) = T_pX$ at every $p \in W_x^u(f) \cap W_y^s(f)$. With the
    above decomposition, we obtain
    \begin{eqnarray*}
        T_pR &=& (T_p X)^+ \\ &=& (T_pW_x^u(f))^+ +
        (T_pW_y^s(f))^+\\
        &=& T_p (W_x^u(f\vert_R)) + T_p(W_y^s(f\vert_R)),
    \end{eqnarray*}
    i.e. $W_x^u(f\vert_R)$ and $W_y^s(f\vert_R)$, viewed as
    submanifolds of $R$, intersect transversely at $p$. Hence
    $(f\vert_R,g\vert_R)$ is Morse-Smale, too.
\end{proof}

\subsection{The $\mathbb Z_2$--homology ring of $R$}
\label{sec:homology-of-R}
We now fix a function $f_\xi: X \to \mathbb R$ and a $\tau$--invariant
Riemannian metric $g$ such that $(f_\xi,g)$ and hence also
$(f_\xi\vert_R,g \vert_R)$ are Morse--Smale. In this subsection we
prove:
\begin{prop}
    \label{prop:Isomorphism-H(X)-H(R)}
    The map $\mathrm{Crit}\,f_\xi \vert_R \to \mathrm{Crit}\,f_\xi$,
    $x \mapsto x$, induces a canonical ring isomorphism
    \begin{equation*}
        H(R; \mathbb Z_2) \to H(X; \mathbb Z_2)
    \end{equation*}
    which doubles the degree, i.e. it takes $H_*(R;\mathbb Z_2)$ to
    $H_{2*}(X;\mathbb Z_2)$.
\end{prop}

This proposition already lets us infer some information about the
quantum homology of $R$. Since $H(X;\mathbb Z_2)$ is generated by
$H_{2n-2}(X;\mathbb Z_2)$ as a ring, it shows that $H_*(R;\mathbb
Z_2)$ is generated by $H_{n-1}(R;\mathbb Z_2)$. Under the assumption
$N_R \geq 2$ on the minimal Maslov number (which we need to assume
anyway), this together with Proposition 6.1.1 in
\cite{biran-cornea_QuantumStrForLag07} implies:
\begin{cor}
    The $\Lambda_R$--wide/narrow dichotomy holds for $R$,
    i.e. $QH(R;\Lambda_R)$ either vanishes or is
    isomorphic to $H(R;\mathbb Z_2) \otimes \Lambda_R$
    as a $\Lambda_R$--module.
\end{cor}

The proof of Proposition \ref{prop:Isomorphism-H(X)-H(R)} will occupy
the rest of this section.  We first concentrate on proving that the
described map induces an isomorphism of $\mathbb Z_2$--vector
spaces. The main ingredient is the following statement.
\begin{lem}
    \label{lemma:total_homology_of_R}
    The rank of the total homology $H(R;\mathbb Z_2)$ is equal to
    the number of critical points of $f_\xi\vert_R$, i.e. $f_\xi\vert_R$ is a
    perfect Morse function.
\end{lem}
This is a special case of a more general result due to Duistermaat,
Theorem 3.1 in \cite{duistermaat--convexity_and_tightness__1983}. We
include the proof for the convenience of the reader, following
\cite{duistermaat--convexity_and_tightness__1983} closely.

\begin{proof}[Proof of Lemma \ref{lemma:total_homology_of_R}]
    Assume without loss of generality that $\xi \in \mathfrak q$ is
    chosen such that $t \mapsto \exp \,t \xi$ is periodic, which is
    the case for a dense subset of $\mathfrak q$. Furthermore normalize
    $\xi$ such that $\exp \,\xi = 1 \in Q$.

    From the relation between $\tau$ and the torus action, $\tau \circ
    q = q^{-1} \circ \tau$ for $q \in Q$, we obtain
    \begin{equation*}
        \tau (\exp (t\xi)
        \cdot x) = \exp (- t \xi) \cdot x
    \end{equation*}
    for $x \in R = \mathrm{Fix}\,\tau$. It follows that the
    ``half-turn map'' $\exp (\frac{1}{2} \xi) = \exp (- \frac{1}{2}
    \xi)$ defines an involution of $R$, whose fixed point set we
    denote by $R_{(1)}$. It is clear that $R_{(1)}$ contains the
    critical set of $f_\xi \vert_R$ (the fixed point set of the
    $Q$--action), but it might be strictly bigger. Therefore, we set
    more generally
    \begin{equation*}
        R_{(k)} := \big\{x \in R_{(k-1)} ~\big \vert~ \exp
        (2^{-k} \xi ) \cdot x = x\big\}
    \end{equation*}
    for $k \geq 1$. This defines a decreasing sequence of submanifolds
    of $R$ which eventually stabilizes. Since $x \in R_{(\ell)}$ for
    all $\ell$ implies $d f_\xi(x) = 0$, we see that in fact we have
    $R_{(\ell)} = \mathrm{Crit}\, f_\xi \vert_R$ for all sufficiently
    large $\ell$.

    Now each $R_{(\ell)}$ is the fixed point set of an involution of
    $R_{(\ell-1)}$. It was shown by Floyd
    \cite{Floyd--On_periodic_maps_1952} as an application of Smith
    theory (see e.g. Section III.3 in
    \cite{Bredon_Intro-Compact-Transf-Groups_1972}) that if $M$ is a
    compact manifold and $N$ the fixed point set of a periodic map $M
    \to M$ of prime period $p$, then $\mathrm{dim}\, H(M;\mathbb Z_p)
    \geq \mathrm{dim}\, H(N;\mathbb Z_p)$ (the inequality does not
    necessarily hold in individual degrees). As a consequence, we
    obtain
    \begin{equation*}
        \mathrm{dim}\,H(R;\mathbb Z_2) \geq 
        \mathrm{dim}\,H(R_{(1)}; \mathbb Z_2) \geq \dots \geq
        \mathrm{dim}\,H(\mathrm{Crit}\,f_\xi \vert_R; \mathbb Z_2).
    \end{equation*}
    Since the function $f_\xi \vert_R$ is Morse, the right-hand side
    is simply the number of its critical points. But this is $\geq
    \mathrm{dim}\,H(R; \mathbb Z_2)$, which concludes the proof.
\end{proof}

Since we know that $f_\xi$ is a perfect Morse function for $X$, this
lemma together with Lemma \ref{lemma:Critical_points_on_X_and_R} allows
to conclude that the map described in Proposition
\ref{prop:Isomorphism-H(X)-H(R)} induces an isomorphism $H(R;\mathbb
Z_2) \cong H(X;\mathbb Z_2)$ of vector spaces.

We next argue that the isomorphism preserves the ring structure,
i.e. is compatible with the intersection product on both sides. Recall
first that in order to define the intersection product Morse
homologically, one has to work with two different Morse functions $f,
f'$ and a metric $g$ satisfying a version of Morse--Smale
transversality. The product is then induced by the map
\begin{equation*}
    \mathbb Z_2 \langle \mathrm{Crit}\, f \rangle \otimes \mathbb Z_2
    \langle\mathrm{Crit}\,f' \rangle \to \mathbb Z_2 \langle
    \mathrm{Crit} f \rangle,
    \quad x \otimes y \mapsto \sum
    \langle x \ast y,z\rangle z,
\end{equation*}
where the sum runs over all $z \in \mathrm{Crit}\,f'$ such that $\vert
x \vert + \vert y \vert - \vert z \vert = n$ and the coefficient
$\langle x \ast y,z\rangle \in \mathbb Z_2$ counts element of the
moduli space $\mathcal Y(x,y,z;f,f')$, whose elements are
configurations consisting of a flow line of $-\nabla f$ from $x$ to
$z$ which is hit by a flow line of $-\nabla f'$ emanating from $y$.

In our particular case, we work again with $\tau$--invariant Morse
functions having all their critical points on $R$. While it is
tempting to choose $f = f_\xi$ and $f' = f_{\xi'}$ for certain $\xi,
\xi' \in \mathfrak q$, this would be problematic because the critical
sets of these coincide, which obstructs the general position
requirement. An easy solution to this is to choose instead a small
perturbation $f' = \tilde f_{\xi'}$ obtained from $f_{\xi'}$ by
pushing the critical points slightly off their original position via a
$\tau$--equivariant diffeomorphism supported on a small neighborhood
of the critical set. An argument similar to the (omitted) proof of
Proposition \ref{prop:tau-invt-Morse-Smale} shows that Morse-Smale
transversality holds for an ``even more generic'' $\tau$-invariant
metric $g$. With these choices, $\tau$ induces an involution
\begin{equation*}
    \tau_* : \mathcal Y(x,y,z;f_\xi,\tilde f_{\xi'}) \to \mathcal
    Y(x,y,z;f_\xi,\tilde f_{\xi'})
\end{equation*}
on the moduli spaces used to compute the intersection product on
$H(X;\mathbb Z_2)$, and the fixed points of $\tau_*$ are precisely
those configurations which are contained in $R$. All configurations
which are \emph{not} entirely contained in $R$ occur in pairs, so that
\begin{equation*}
    \#_{\mathbb Z_2} \mathcal Y(x,y,z;f_\xi,\tilde f_{\xi'}) \,=\, \#_{\mathbb
        Z_2} \mathcal Y(x,y,z;f_\xi\vert_R,\tilde f_{\xi'}\vert_R).
\end{equation*}
Because the moduli space on the right-hand side contains precisely
those configurations counted for the intersection product on
$H(R;\mathbb Z_2)$, this finishes the proof of the statement about the
ring structure and hence of Proposition
\ref{prop:Isomorphism-H(X)-H(R)}.  \hfill $\square$

\section{Real Discs and their Doubles}
\label{sec:real-discs-real}

In this section we study the holomorphic curves that enter in the
definition of the quantum invariants of $R$ and $X$, where holomorphic
is always understood with respect to the standard integrable complex
structure $J_0$.

\begin{defn}
    A \emph{real disc} is a holomorphic map $u: (D^2,\partial D^2) \to
    (X,R)$, where $D^2 \subset \mathbb C $ denotes the closed unit
    disc. A \emph{real rational curve} is a holomorphic map $v:
    \mathbb C P^1 \to X$ such that $\tau \circ v \circ c = v$, with
    $c: \mathbb C P^1 \to \mathbb C P^1$ denoting complex
    conjugation.
\end{defn}

The $\tau$--anti-invariance of $J_0$ implies that for every real disc
$u$, the disc
\begin{equation*}
    \overline{u} : (D^2,\partial D^2) \to (X,R),\quad z \mapsto \tau (u
    (\overline z)),
\end{equation*}
obtained by Schwarz reflection is again holomorphic. Gluing the discs
$u$ and $\overline u$ along their common boundary yields a holomorphic sphere
\begin{equation*}
    u^\sharp  : \mathbb C P^1 \to X
\end{equation*}
which we refer to as the \emph{double} of $u$. In the usual
identification of $D^2$ with the extended upper half-plane $\mathbb H
\cup \{\infty\} \subset \mathbb C P^1$, the double is given by
\begin{equation*}
    u^\sharp (z) \,=\,
    \begin{cases}
        u(z),& z \in \mathbb H \cup \{\infty\},\\
        \tau(u(\overline z)),& \text{otherwise}.
    \end{cases}
\end{equation*}

The main point of this section is to provide a representation of real
discs and their doubles in toric homogeneous coordinates which will be
used for proving the regularity of $J_0$ for these curves in Section
\ref{sec:transversality}. This strategy is very similar to what Cho
and Oh did for discs with boundary on a torus fibre and was actually
inspired by their work \cite{Cho-Oh_FloerCohomology2006}.

\subsection{Lifting real discs and their doubles}
\label{sec:homog-coord-form}
Recall that we constructed the toric manifold $X$ as a quotient of an
open set $\mathcal U \subset \mathbb C^N$ by some torus action. It is
not difficult to see that any real disc $u$ can be lifted
holomorphically to $\mathcal U$, i.e. that there exists a holomorphic
map $w$ such that the diagram
\begin{diagram}
    &&(\mathbb C^N, \pi^{-1}(R))\\
    &\ruTo^{w} & \dTo{\pi}\\
    (D^2, \partial D^2)& \rTo{\quad u} & (X,R)
\end{diagram}
commutes. For the purpose of our transversality proof we would ideally
like such a lift to take boundary values in some totally real
submanifold of $\mathbb C^N$ which projects to $R$ under $\pi:
\mathcal U \to X$ and whose tangent bundle we understand. In our case
the obvious candidate for this submanifold would be $\mathbb R^N$, but
for a non-constant $u$ it is not possible to find a holomorphic lift
satisfying $w(\partial D^2) \subset \mathbb R^N$ by Liouville's
theorem. In this respect our situation is somewhat different from that
considered by Cho and Oh in \cite{Cho-Oh_FloerCohomology2006}, who
could lift discs with boundary on a torus fibre such that the boundary
was mapped to a totally real torus in $\mathbb C^N$.

Our solution to this problem is to take out a boundary point and lift
only the restriction of $u$ to $\mathbb H \subset D^2$, or
equivalently the restriction of the double $u^\sharp$ to $\mathbb C
\subset \mathbb C P^1$. It will turn out that this is sufficient
for a transversality proof similar to Cho and Oh's to work.

\begin{prop}
    \label{thm::hom_coords}
    Let $u: (D^2,\partial D^2) \to (X,R)$ be a real disc and let
    $u^\sharp: \mathbb C P^1 \to X$ be its double. The restriction of
    $u^\sharp$ to $\mathbb C \subset \mathbb C P^1$ admits a lift
    \begin{equation*}
        w =(w_1,\dots,w_N): \,(\mathbb C, \mathbb R)  \,\to \, (\mathbb C^N, \mathbb R^N),
    \end{equation*}
    where the $w_i: \mathbb C \to \mathbb C $ are holomorphic
    functions of the form
    \begin{equation*}
        \textstyle
        w_i(z) = a_i \cdot \prod_{j=1}^{\alpha_i} (z - p_{i,j}) (z - \overline{p_{i,j}}) \cdot
        \prod_{k = 1}^{\beta_i} (z - q_{i,j}),
    \end{equation*}
    with constants $a_i \in \mathbb R $, integers $\alpha_i, \beta_i \geq 0$,
    and points $p_{i,j} \in \mathbb C \smallsetminus \mathbb R $,
    $q_{i,j} \in \mathbb R $.
\end{prop}

\begin{proof}
    We start with an arbitrary holomorphic lift
    \begin{equation*}
        w = (w_1,\dots,w_N): (\mathbb C, \mathbb R) \to (\mathbb C^N, \pi^{-1}(R))
    \end{equation*}
    of the double $u^\sharp: \mathbb C \to X$, which exists because
    the pullback of the holomorphic principal $K_{\mathbb C}$--bundle
    $\mathcal U \to X$ to $\mathbb C $ by $u^\sharp$ is
    holomorphically trivial (see
    e.g. \cite{Grauert--Analytische_Faserungen_1958}).  Denote by $I_0
    \subset [1,N]$ the set of all $i$ such that $u (D^2)$ is contained
    in the hypersurface $D_i$, and set $\overline{I_0} := [1,N] \smallsetminus
    I_0$. For $i \in \overline{I_0}$, $u^{\sharp}$ hits $D_i$ at only finitely
    many points, and hence $w_i$ has only finitely many
    zeros. Clearly, the set of zeros of $w_i$ is invariant under
    complex conjugation by the $\tau$--invariance of
    $u^{\sharp}$. Denote by $p_{i,1}, \overline{p_{i,1}}, \dots,
    p_{i,\alpha_i}, \overline{p_{i,\alpha_i}}$ the zeros of $w_i$ in
    $\mathbb C \smallsetminus \mathbb R$, and by $q_{i,1},\dots,
    q_{1,\beta_i}$ those in $\mathbb R $. In both cases, a zero with
    multiplicity $m$ is treated as $m$ zeros.

    Now define a new holomorphic map
    \begin{equation*}
        w' =   (w_1',\dots,w_N') : \mathbb C \to \mathbb C^N
    \end{equation*}
    by setting
    \begin{equation*}
        w_i'(z) \, := \, \frac{w_i(z)}{\chi_i(z)}
    \end{equation*}
    with $\chi_i(z) \,:=\, \prod_{j=1}^{\alpha_i} (z- p_{i,j}) (z -
    \overline {p_{i,j}}) \cdot \prod_{j=1}^{\beta_i} (z - q_{i,j})$.
    Finally set
    \begin{equation*}
        u' \,:=\, \pi \circ w': \mathbb C \to X.
    \end{equation*}

    \begin{claim*}
        $u': \mathbb C \to X$ is constant with image in $R$.
    \end{claim*}
    
    It is easy to see that $u'$ extends over $\infty$ to a holomorphic
    map $u': \mathbb C P^1 \to X$, using e.g. the projective embedding
    $\Psi: X \to \mathbb C P^L$ described in Section
    \ref{sec:sympl-toric-manif}.
    Notice that $\infty$ is the only point where this extension can
    possibly meet any of the hypersurfaces $D_i$, $i \in
    \overline{I_0}$, because for those $i$ the corresponding $w_i' :
    \mathbb C \to \mathbb C$ are nowhere-vanishing. It follows that
    all the $D_i$ with $u'(\mathbb C P^1) \cap D_i \neq \emptyset$,
    including those with $i \in I_0$, intersect at $u'(\infty)$. In
    particular, this intersection is non-empty, and thus the cone
    $\sigma$ generated by the corresponding normal vectors $v_i$ of
    $\Delta$ is contained in the normal fan of $\Delta$. Let $\sigma'
    \in \Sigma^{(n)}$ be an $n$--dimensional cone containing $\sigma$
    as a face. Hence $u'(\mathbb C P^1)$ is entirely
    contained in the set $V_{\sigma'}$, which by Lemma
    \ref{lem::affine_open_sets} is biholomorphic to $\mathbb C^n$. So
    $u'$ is constant by Liouville's theorem, and the fact that all
    $\chi_j$ send real numbers to real numbers implies that the image
    of $u'$ actually lies in $R$. This proves the claim.

    To conclude the proof of the proposition, note that we may assume
    that also the lift $w'$ of $u'$ is constant (otherwise $w'$ is of
    the form $w'(z) = w''(z) \cdot w'(0)$ with $w''$ taking values in
    $K_{\mathbb C}$, in which case we can replace the original $w$ by
    $(w'')^{-1} \cdot w$), and consequently the image of $w'$ is a point in
    $\pi^{-1}(R)$. By eventually acting with some element of
    $K_{\mathbb C}$, we may further assume that this point actually
    lies in $\mathbb R^N$. Recalling the way $w'$ and $w$ are related
    shows that $w$ is of the claimed form.
\end{proof}

\subsection{Reading off the Maslov index}
\label{sec:reading-maslov-index}
As for the relationship between the Maslov index $\mu(u)$ of a
real disc and the first Chern class $c_1(u^\sharp)$ of its double, we have
\begin{eqnarray*}
    c_1(u^\sharp) &=& \frac{1}{2} \big( \mu(u) + \mu(\overline u) \big) \\
    &=& \mu(u).
\end{eqnarray*}
The first equality is well-known
, and the second follows from $\mu (u) = \mu(\overline u)$. Given that
$c_1(X)$ is Poincar\'e dual to $[D_1] + \dots + [D_N]$, we conclude
that the Maslov index of $u$ is equal to the intersection number of
$u^\sharp$ with the toric divisors. In case $u$ is not entirely
contained in any $D_i$ and hence intersects each of them only a finite
number of times, this number is simply
\begin{equation*}
    \mu(u) = \sum_{i = 1}^N \# \big\{\text{zeros of $w_i$}\big\}
\end{equation*}
where the $w_i$ are the components of any holomorphic lift $w: \mathbb
C \to \mathbb C^N$ of $u^\sharp$, counting with multiplicities. We
assume here that $u(\infty)$ is not contained in any $D_i$, so that
all intersection points of $u$ with the $D_i$ correspond to zeros of
the $w_i$.

If $u$ is entirely contained in some $D_i$, or equivalently $w_i$
vanishes identically, one has to represent $[D_i]$ by a linear
combination of those $[D_j]$ for which $u$ is not contained in $D_j$
in order to determine the intersection number of $u^\sharp$ with
$[D_i]$.

This results in the following general recipe for writing the Maslov
index in terms of a weighted count of the numbers of isolated zeros of
the components of $w$. As in the proof of Proposition
\ref{thm::hom_coords} we denote by $I_0 = \{i_1,\dots,i_\ell\}$ the
set of those $i \in [1,N]$ such that $u(D^2)$ is contained in $D_i$ and
set $\overline{I_0} := [1,N] \smallsetminus I_0$. Moreover, we assume
that $u(\infty) \notin \bigcup_{i \in \overline {I_0}} D_i$, which can
be arranged using a suitable reparametrisation. Observe now that
$\bigcap_{i \in I_0} D_i \neq \emptyset$ because $u(D^2)$ is contained
in the intersection. It follows that the $ v_{i_1},\dots,v_{i_\ell}$
generate a cone in $\Sigma$, and so they can be extended to a
$\mathbb Z $--basis of $\mathbb Z^n$. Fixing any such basis, we denote
by $\varepsilon_{i_1},\dots,\varepsilon_{i_\ell} \in (\mathbb Z^n)^*$
those elements of the corresponding dual basis of $(\mathbb Z^n)^*$
which are dual to the $v_{i_1},\dots,v_{i_\ell}$, i.e. $\langle
\varepsilon_{i_j}, v_{i_k} \rangle = \delta_{jk}$ for $j,k =
1,\dots,\ell$.

\begin{prop}
\label{prop:Maslov-index-formula}
 The Maslov index of $u$ can be expressed as 
    \begin{equation*}
        \mu(u) = \sum_{j \in \overline{I_0}} \Big(1 - \sum_{i \in I_0} \langle
            \varepsilon_i, v_j \rangle \Big) (2 \alpha_j + \beta_j),
    \end{equation*}
    where the $\alpha_i, \beta_i$ are as in Proposition
    \ref{thm::hom_coords}.
\end{prop}

\begin{proof} By the description of $H(X;\mathbb Z_2)$ as given in
    Section \ref{sec:sympl-toric-manif}, the identity
    \begin{equation*}
        \sum_{j = 1}^N \langle \varepsilon_i, v_j \rangle [D_j] = 0
    \end{equation*}
    in $H_{2n-2} (X; \mathbb Z_2)$ holds for all $i \in I_0$. The
    choice of the $\varepsilon_{i}$ implies that
    \begin{equation*}
        [D_i] = - \sum_{j \in \overline{I_0}} \langle \varepsilon_i, v_j \rangle [D_j]
    \end{equation*}
    for all $i \in I_0$. Since $c_1(X)$ is Poincar\'e dual to $[D_1] +
    \dots +[D_N]$, we obtain
    \begin{eqnarray*}
        c_1(X) &=& PD \left(\sum_{i \in I_0} \Big(-\sum_{j \in \overline{I_0}}
        \langle \varepsilon_i,v_j \rangle [D_j] \Big) + \sum_{j \in
            \overline{I_0}} [D_j]\right)\\
            &=& PD \left( \sum_{j \in\overline{ I_0}} \Big(1 - \sum_{i \in I_0}
            \langle \varepsilon_i, v_j \rangle \Big) [D_j] \right),
    \end{eqnarray*}
    which implies the formula claimed for $\mu(u)$.
\end{proof}

\begin{example*}
    We include a baby example to illustrate how reading off the Maslov
    index works. Consider the fan $\Sigma$ in $\mathbb R^2$ whose
    1-dimensional cones are generated by
    \begin{equation*}
        v_1 = (1,0), ~v_2 = (0,1),~ v_3 = (-1,-1), ~v_4 = (0,-1).
    \end{equation*}
    The corresponding toric manifold is $X = \mathbb C P^2 \#
    \overline{\mathbb C P^2}$, a toric one-point blow-up of $\mathbb C
    P^2$. Observe that we have $[z_1:z_2:z_3:z_4] = [tz_1,stz_2,tz_3,sz_4]$ for
    $(s,t) \in T^2_{\mathbb C}$ because the map $\pi: \mathbb Z^4 \to
    \mathbb Z^2$, $e_i \mapsto v_i$, has kernel $\mathbb K = \mathbb Z
    \big\langle (0,1,0,1),$ $(1,1,1,0) \big\rangle$. Consider now the
    real rational curve $u$ whose restriction to $\mathbb H$ is given
    by
    \begin{equation*}
        u(z) = \big[z:1:1:0\big],
    \end{equation*}
    or equivalently, whose lift to $\mathbb C^2$ is $w(z) =
    (z,1,1,0)$.  It sweeps out one half of the divisor $D_4$. For $z
    \neq 0$ we can write $u(z) = [1:1:z^{-1}:0]$ by acting with
    $(1,z^{-1}) \in T^2_{\mathbb C}$, and hence $u$ intersects $D_3$
    at $\infty$. To move this intersection away from $\infty$, we
    reparametrise $u$ using the element $\phi \in \mathrm{Aut}(D^2)$,
    $z \mapsto \frac{z}{z-1}$, such as to obtain the new real disc
    \begin{equation*}\textstyle
        u'(z) = \big[\frac{z}{z-1} : 1 : 1 : 0 \big] = \big[z : 1 :
        z-1:0 \big],
    \end{equation*}
    which has the form required for the formula of Proposition
    \ref{prop:Maslov-index-formula} to be applicable. We have $I_0 =
    \{4\}$, $\alpha_i = 0$ for all $i$, $\beta_1 = \beta_3 = 1$ and
    $\beta_2 = 0$. Choosing e.g. the vector $(1,1)$ to extend $v_4 =
    (0,-1)$ to a basis of $\mathbb Z^2$, the required element of the
    dual basis is $\varepsilon_4 = (1,-1)$, where we identify
    $(\mathbb Z^2)^*$ with $\mathbb Z^2$ in the usual way. With
    $\langle \varepsilon_4,v_1\rangle = 1$ and $\langle
    \varepsilon_4,v_3 \rangle = 0$ we obtain
    \begin{equation*}
        \mu(u) = \mu(u') = (1 - \langle \varepsilon_4,v_1\rangle)\beta_1 + (1
        - \langle \varepsilon_4,v_3\rangle )\beta_3 = 1
    \end{equation*}
    as expected, because $D_4$ is the exceptional divisor.\footnote{In
        particular we have $C_X = 1$, so that our main theorem is not
        applicable to $X$.}
\end{example*}

\section{Regularity of the Standard Complex Structure}
\label{sec:transversality}

In this section we will prove the following theorem.
\begin{thm}
    \label{thm::transversality}
    The standard complex structure $J_0$ is regular for all
    holomorphic real discs $u: (D^2,\partial D^2) \to (X,R)$ and
    all holomorphic spheres $v: \mathbb C P^1 \to X$.
\end{thm}

The proof is similar to Cho and Oh's regularity proof in
\cite{Cho-Oh_FloerCohomology2006}. We will give full details only for
the statement about discs, and then sketch the essentially identical
proof for the sphere case in the last subsection.

\subsection{Reformulating the problem}
Recall that Fredholm regularity of $J_0$ for $u$ means that the
linearisation $D_u \overline \partial_{J_0}$ of the operator
$\overline \partial_{J_0}$ at $u$ is surjective. Here we view
$\overline \partial_{J_0}$ as a section in a suitable Banach manifold
setting. In our case, this linearisation is an operator
\begin{equation*}
    D_u \overline \partial_{J_0}: W^{1,p} (u^* TX, u^*TR) \,\to\, L^p(\Lambda^{0,1}T^*D\otimes u^*TX)
\end{equation*}
for some $p > 2$, whose restriction to smooth\footnote{Here
    and in the following ``smooth'' means ``smooth in the interior and
    all derivatives continuous up to the boundary''.} sections
is the standard Dolbeault operator
\begin{equation*}
    \overline \partial: C^{\infty} (u^*TX,u^*TR) \,\to\,
    \Omega^{0,1} (u^*TX).
\end{equation*}
Since we know that $D_u \overline \partial_{J_0}$ is Fredholm by
general theory and because smooth sections are dense, it suffices to
show that $\overline \partial$ is surjective, i.e.
\begin{equation*}
    \mathrm{coker}\,\overline \partial = 0.
\end{equation*}

This is equivalent to the vanishing of a certain sheaf cohomology
group. To first discuss this in general terms, let $(E,F)$ be a
holomorphic \emph{bundle pair} over $(D^2,\partial D^2)$, that is, let
$E$ be a holomorphic vector bundle over $D^2$ and let $F$ be a totally
real subbundle over $\partial D^2$. Denote by $\mathcal O (E,F)$ the
sheaf of holomorphic sections of $E$ taking boundary values in $F$, by
$\mathcal A^0(E,F)$ the sheaf of smooth sections of $E$ taking
boundary values in $F$, and by $\mathcal A^{(0,1)}(E)$ the sheaf of
smooth $E$--valued $(0,1)$--forms. The latter two are fine sheaves,
i.e. they admit partitions of unity.

\begin{prop}
    \label{prop::fine-resolution}
    The sequence of sheaves $\mathcal A^0(E,F)
    \xrightarrow[]{\overline \partial} \mathcal A^{(0,1)}(E) \to 0 \to
    \dots$ is a fine resolution of $\mathcal O(E,F)$.
\end{prop}

This is proved (for arbitrary Riemann surfaces with boundary) e.g. in
\cite{Katz-Liu--Enumerative_geometry__2001}. It follows that the sheaf
cohomology groups of $\mathcal O(E,F)$, which are by definition the
right derived functors of the global section functor applied to
$\mathcal O(E,F)$, satisfy
\begin{equation*}
    H^*(D^2;\mathcal O(E,F)) ~=~
    \begin{cases}
        \mathrm{ker}\, \overline \partial,& \quad * = 0\\
        \mathrm{coker}\, \overline \partial,& \quad * = 1\\
        0,& \quad * \geq 2.
    \end{cases}    
\end{equation*}

The idea for our specific problem is to show $H^1(D;\mathcal O(E,F)) =
0$ for $(E,F) = (u^*TX,u^*TR)$ by exhibiting a suitable acyclic cover
of $D^2$ and making use of the fact that sheaf cohomology groups can
be computed as sheaf-valued \v{C}ech cohomology groups with respect to
an acyclic cover.

\subsection{An acyclic cover}
\label{sec:finding-an-acyclic} 
By definition, \emph{acyclicity} of a cover for a sheaf means that all
higher cohomology groups ($H^i$ with $i \geq 1$) of the sheaf
restricted to arbitrary intersections of the covering sets
vanish. Consider the cover $\mathfrak U = \{U_0,U_1\}$ of $D^2$ with
\begin{equation*}
    U_i \,=\, D^2 \smallsetminus \{z_i\},
\end{equation*}
where $z_0 = 0$ and $z_1 = \infty$, thinking of $D^2$ as the extended
upper half-plane. The purpose of this subsection is to prove that
$\mathfrak U$ is acyclic for $\mathcal O(E,F)$, i.e. that
\begin{equation*}
    H^1(U_0; \mathcal O(E,F)) = H^1(U_1; \mathcal O(E,F)) = H^1(U_0
    \cap U_1; \mathcal O(E,F)) = 0.
\end{equation*}
As both $U_i$ are biholomorphic to the upper half-plane, in the
following we simply write $\mathbb H$ to represent either of
them. Moreover, we write $\mathbb H^*$ for $U_0 \cap U_1$.

Let $w: (\mathbb H, \mathbb R) \to (\mathbb C^N, \mathbb R^N)$ be a
lift of $u\vert_{\mathbb H}$, which exists by Proposition
\ref{thm::hom_coords}, and consider the following bundle pairs over
$(\mathbb H, \mathbb R)$:
\begin{eqnarray*}
    (E_K, F_K) &:=& (w^*(T Orb_{K_{\mathbb C}}), w^*(T Orb_{K_{\mathbb C}} \cap \mathbb R^N)),\\
    (\tilde E, \tilde F) &:=& (w^*\mathbb C^N, w^* \mathbb R^N),\\
    (E,F) &:=& (u^*TX, u^*TR).
\end{eqnarray*}
Here $T Orb_{K_{\mathbb C}}$ denotes the subbundle of $T\, \mathbb
C^N$ whose fibre at $z \in \mathbb C^N$ is the tangent space of the
$K_{\mathbb C}$--orbit through $z$. The following statement is proven
in exactly the same way as Lemma 6.3 in
\cite{Cho-Oh_FloerCohomology2006}.
\begin{lem}
    \label{lem:exact-sheaf-seq}
    The natural sequence of sheaves
    \begin{equation*}
        0 \to \mathcal O (E_K,F_K) \to \mathcal O (\tilde E,
        \tilde F) \to \mathcal O( E,F) \to 0
    \end{equation*}
    is exact.
\end{lem}

This implies that there exists a long exact sequence in cohomology of the
form
\begin{align*}
    0 &\to H^0(\mathbb H; \mathcal O (E_K,F_K)) \to
    H^0(\mathbb H; \mathcal O (\tilde E,\tilde F)) \to H^0(\mathbb H;\mathcal O (E,F))\\
    &\to H^1(\mathbb H;\mathcal O (E_K,F_K)) \to H^1(\mathbb
    H;\mathcal O (\tilde E, \tilde F)) \to H^1(\mathbb H;\mathcal O
    (E,F)) \to 0.
\end{align*}
The sequence stays exact when $\mathbb H$ is replaced by
the open subset $\mathbb H^*$. To show $H^1(\mathbb H; \mathcal O( E,
F)) = H^1(\mathbb H^*;\mathcal O( E, F)) = 0$ and thus
the acyclicity of $\mathfrak U$, it is hence sufficient to show
\begin{equation*}
    H^1(\mathbb H; \mathcal O(\tilde{E}, \tilde{F})) = 
    H^1(\mathbb H^*; \mathcal O(\tilde{ E}, \tilde{F})) = 0.
\end{equation*}
Since $(\tilde E, \tilde F) \,=\, (\mathbb H \times \mathbb C^N,
\mathbb R \times \mathbb R^N)$ is the trivial bundle with a trivial
totally real subbundle, this follows from the next proposition.
\begin{prop}
    \label{prop::surjectivity_trivial_bundle}
    Given a smooth 
    function $g: \mathbb H \to \mathbb C$ such that $g(\mathbb R )
    \subset \mathbb R $, there exists a smooth function $f: \mathbb H
    \to \mathbb C $ satisfying $\overline \partial f = g$ and
    $f(\mathbb R ) \subset \mathbb R $.
\end{prop}
We include (parts of) the proof for the convenience of the reader,
essentially reproducing the proof of Theorem 1.11 in
\cite{Gunning--Function_theory_on_RS}.

\begin{proof}
    For compactly supported functions $g$ the statement is a rather
    well-known result from complex analysis.

    Suppose now that $g$ is not compactly supported. For $n \in
    \mathbb N$, let $B_n := \big\{z \in \mathbb H ~\big\vert~ \vert z
    \vert < n\big\}$ be the open half-disc of radius $n$, let $\rho_n$
    be a cut-off function with $\rho_n \equiv 1$ on $\overline {B_n}$,
    $\mathrm{supp}\,\rho_n \subset B_{n+1}$ and let $g_n := \rho_n
    \cdot g$. We will construct a sequence of smooth functions $f_n:
    \mathbb H \to \mathbb C$ such that
    \renewcommand{\theenumi}{\roman{enumi}}
    \begin{enumerate}
    \item $\overline \partial f_n = g$ in $B_n$,
    \item $f_n (\mathbb R ) \subset \mathbb R $,
    \item $\vert f_{n+1} (z) - f_{n}(z) \vert \leq 2^{-n}$ for all $z \in
        \overline {B_{n-1}}$.
    \end{enumerate}
    Suppose that functions $f_1,\dots,f_n$ satisfying these conditions
    have been constructed. To define $f_{n+1}$, start with a function
    $\varphi: \mathbb H \to \mathbb C$ solving $\overline \partial
    \varphi = g_{n+1}$ and $\varphi (\mathbb R)\subset \mathbb R$,
    which exists because $g_{n+1}$ has compact support. The function
    $\varphi - f_n$ is holomorphic on $B_n$ and satisfies $(\varphi -
    f_n)(\mathbb R) \subset \mathbb R $; by truncating its power
    series expansion around 0 at a sufficiently high power, we obtain
    a polynomial $\psi$ which is holomorphic on $B_{n+1}$, satisfies
    $\psi(\mathbb R ) \subset \mathbb R $, and $\big\vert \varphi(z) -
    f_n(z) - \psi (z) \big\vert < 2^{-n}$ for all $z \in
    \overline{B_{n-1}}$. Then we define
    \begin{equation*}
        f_{n+1} \,:= \, \varphi - \psi
    \end{equation*}
    and note that it has the required properties.

    As for the convergence of this sequence of functions, observe that
    all $f_m$ with $m \geq n$ satisfy $\overline \partial f_m = g$ on
    $B_n$. We can hence write
    \begin{equation*}
        f_m(z) = f_n(z) + \sum_{j=n+1}^{m} f_{j}(z) - f_{j-1}(z)
    \end{equation*}
    for $z \in B_n$. Since all $f_j - f_{j-1}$ in the sum are
    holomorphic on $B_n$ and because of the uniform bounds they
    satisfy on $\overline{B_{n-1}}$, the sum converges for $m \to
    \infty$ to a function which is holomorphic on $B_n$ and takes
    $\mathbb R $ to $\mathbb R $. Since this is true for every $n$, it
    follows that the $f_m$ converge for $m \to \infty$, uniformly on
    compact sets, to a function $f: \mathbb H \to \mathbb C $
    satisfying $\overline \partial f = g$ and $f(\mathbb R ) \subset
    \mathbb R $.
\end{proof}

\subsection{A \v{C}ech cohomology computation}
\label{sec:cech-cohom-comp}
In this section we finish the proof of
\begin{equation*}
    H^1(D^2; \mathcal O(E, F)) =  0
\end{equation*}
and hence of Theorem \ref{thm::transversality}. We make use of the
following classical theorem:
\begin{thm}
    [Leray] Let $X$ be a topological space, $\mathcal S$ a sheaf on
    $X$, and let $\mathfrak U = \{U_i\}_{i \in I}$ be a countable cover of $X$
    which is acyclic with respect to $\mathcal S$
    . Then
    \begin{equation*}
        H^*(X;\mathcal S) \,\cong\, \check{H}^*(\mathfrak U; \mathcal S).
    \end{equation*}
\end{thm}

Here $H^*(X;\mathcal S)$ denotes the sheaf cohomology of $\mathcal S$
in the sense of ``right derived functor of the global section
functor'', and $\check{H}^*(\mathfrak U; \mathcal S)$ denotes \v{C}ech
cohomology with respect to the cover $\mathfrak U$ with values in
$\mathcal S$. So we are done if we can show
\begin{equation*}
    \check{H}^1(\mathfrak U; \mathcal
    O( E,F)) = 0
\end{equation*}
for the acyclic cover $\mathfrak U = \{U_0,U_1\}$ from the previous
subsection. That is, we have to show that every \v{C}ech 1--cocycle is
a coboundary. Since our cover consists of only two sets, a 1--cocycle
is simply a section $\sigma$ of $\mathcal O(E, F)$ over $U_0 \cap U_1
= \mathbb H^*$, and for $\sigma$ to be a coboundary means that there
exist sections $\eta_0, \eta_1$ over $U_0,U_1$ such that
\begin{equation}
    \label{cech-coboundary-cond}
    \sigma \,=\, \eta_0 \big\vert_{U_0 \cap U_1} + \eta_1 \big\vert_{U_0 \cap U_1}.
\end{equation}
In other words, all we have to show is the following claim.

\begin{claim*}
    Any section $\sigma$ of $\mathcal O(E,F)$ over $\mathbb H^*$
    decomposes as 
    \begin{equation*}
        \sigma = \eta_0 + \eta_1
    \end{equation*}
    with sections $\eta_i$ over $\mathbb H^*$ such that $\eta_0$
    extends over $0$ and $\eta_1$ extends over $\infty$.
\end{claim*}

\emph{Proof of the Claim, Step 1.} Fix a lift $w: (\mathbb H, \mathbb R) \to (\mathbb
C^N,\mathbb R^N)$ of the restriction of $u$ to $U_0 = \mathbb H$, and
denote by $w$ also the restriction to $\mathbb H^*$. Let $(E_K,F_K)$,
$(\tilde E, \tilde F)$, $(E,F)$ be the corresponding bundle pairs over
$\mathbb H$, $\mathbb H^*$ as in the last subsection. We claim that
$\sigma$ lifts to a section $\tilde \sigma$ of the trivial
bundle pair $(\tilde E, \tilde F)$ over $\mathbb H^*$. In other words,
we claim the map 
\begin{equation*}
    H^0(\mathbb H^*; \mathcal O(\tilde {E}, \tilde { F})) \,\to
    \, H^0(\mathbb H^*; \mathcal O(E, F)),
\end{equation*}
induced by the sequence of sheaves in Lemma \ref{lem:exact-sheaf-seq}
is surjective, which in view of the long exact sequence is equivalent
to saying that 
\begin{equation*}
    H^1(\mathbb H^*; \mathcal O(E_K, F_K)) = 0.
\end{equation*}
That this is indeed the case follows from the next lemma and
Proposition \ref{prop::surjectivity_trivial_bundle} (which also holds
when $\mathbb H$ is replaced by $\mathbb H^*$).

\begin{lem}
    \label{lem::triviality_of_bundle_pair}
    The bundle pair $(E_K, F_K) = (w^*(T Orb_{K_{\mathbb C}}), w^*(T
    Orb_{K_{\mathbb C}} \cap \mathbb R^N))$ is holomorphically
    trivial.
\end{lem}
\begin{proof}
    Let $z \in \mathcal U \subset \mathbb C^N$. Then we have
    \begin{equation*}
        E_{K,z} \,=\, T_z Orb_{K_{\mathbb C}} \, = \, z \cdot
        T_{\underline{1}}Orb_{K_{\mathbb C}},
    \end{equation*}
    where the right-hand side denotes the subspace of $ T_z \mathbb
    C^N \cong\, \mathbb C^N$ obtained by letting $z$ act on the
    subspace $T_{\underline 1} K_{\mathbb C}$ of $T_{\underline 1}
    \mathbb C^N \cong \mathbb C^N$ in the usual way (i.e. $z \cdot
    (t_1,\dots,t_N) = (z_1t_1,\dots,z_Nt_N)$), and $\underline 1 :=
    (1,\dots,1) \in \mathbb C^N$. This follows directly from the fact
    that $K_{\mathbb C}$ acts freely on $\mathcal U$, which is
    equivalent to saying that the map $ K_{\mathbb C} \ni t \mapsto z
    \cdot t$ is a diffeomorphism onto its image. For $z \in \mathcal U
    \cap \mathbb R^N$ the same argument shows
    \begin{equation*}
        F_{K,z} = (T_z Orb_{K_{\mathbb C}}) \cap \mathbb R^N \, = \, z \cdot
        (T_{\underline{1}}Orb_{K_{\mathbb C}} \cap \mathbb R^N).
    \end{equation*}
    It follows that
    \begin{eqnarray*}
        \mathbb H \times T_{\underline 1} Orb_{K_{\mathbb C}}  &\to &
        E_K,\\
        (z, t) &\mapsto& (z, w(z) \cdot t),
    \end{eqnarray*}
    is a holomorphic trivialization of $E_K$ identifying $F_K$ with a
    constant totally real subbundle.
\end{proof}

\emph{Proof of the Claim, Step 2.} Since $(\tilde E, \tilde F)$ is the
trivial bundle pair $(\mathbb H^* \times \mathbb C^N, \mathbb H^*
\times \mathbb R^N)$, the lift $\tilde \sigma$ can be written as
$\tilde \sigma = (\tilde \sigma_1, \dots, \tilde \sigma_N)$, where
each of the $\tilde \sigma_i: \mathbb H^* \to \mathbb C$ admits a
Laurent series expansion
\begin{equation*}
    \tilde \sigma_i (z) \,=\, \sum_{j= -\infty}^\infty a_{i,j} z^{j}
\end{equation*}
with coefficients $a_{i,j} \in \mathbb R$. We set
\begin{equation*}
    \tilde \eta_{0,i}(z) \,:=\, \sum_{j=0}^\infty a_{i,j} z^j,\qquad
    \tilde \eta_{1,i} (z) \,:=\, \sum_{j=-\infty}^{-1} a_{i,j}z^j
\end{equation*}
for $i = 1,\dots,N$, then $\tilde \eta_0 := (\tilde \eta_{0,1}, \dots,
\tilde \eta_{0,N})$, $\tilde \eta_1 := (\tilde \eta_{1,1}, \dots,
\tilde \eta_{1,N})$, and finally
\begin{equation*}
    \eta_0(z) \,:=\, D \pi\big\vert_{w(z)} \tilde \eta_0(z),\qquad
    \eta_1 (z) \,:=\, D \pi\big\vert_{w(z)} \tilde \eta_1(z).
\end{equation*}
The $\eta_i$ are sections of $(E,F)$ over $\mathbb H^*$ which satisfy
\eqref{cech-coboundary-cond} by construction. We still have to show
that they extend in the required way, which by the removal of
singularities theorem requires only to show that the limits
$\eta_0(z)$ for $z \to 0$ and $\eta_1(z)$ for $z \to \infty$
exist. This is clear for $\eta_0$, because for $z \to 0$
both $\tilde \eta_0(z)$ and $w(z)$ converge, the latter by the way $w$
was chosen. As for $\eta_1$, the problem is that while $\tilde
\eta_1(z)$ converges to some $\tilde \eta_1(\infty)$, $w(z)$ does not
converge for $z \to \infty$. Nevertheless, we are in good shape, as
the following lemma shows.

\begin{lem}
    $D\pi\vert_{w(z)} \tilde \eta_1(z)$ converges to $0 \in T_{u(\infty)}R$ for $z \to
    \infty$.
\end{lem}
\begin{proof}
    We start by constructing a function $\theta: \mathbb H \to
    K_{\mathbb C}$ such that
    \begin{equation*}
        \big\Vert \theta(z) \cdot w(z) \big\Vert \,<\, C
    \end{equation*}
    for some $C > 0$ and all $z \in \mathbb H $. To do so, choose
    first an element $\lambda = (\lambda_1,\dots,\lambda_N) \in
    \mathbb K$ such that all $\lambda_i$ are $> 0$. This
    is possible, since for every $i \in [1,N]$ the vector $-v_i$
    lies in some cone $\sigma \in \Sigma$ by the completeness of
    $\Sigma$ and can hence be expressed as a linear combination of the
    generators $v_{j_1},\dots,v_{j_{m_i}}$ of $\sigma$ with
    \emph{non-negative} coefficients. In other words, we have $v_i +
    \sum_{k=1}^{m_i} a_{j_k} v_{j_k} = 0$ for certain $a_{j_k} \geq
    0$. Setting $\lambda^{(i)} \,:=\, e_i + \sum_{k=1}^{m_i} a_{j_k}
    e_{j_k}$, where the $e_i$ denote the elements of the standard
    basis of $\mathbb Z^N$, and finally $\lambda \,:=\, \sum_{i=1}^N
    \lambda^{(i)}$ yields the required vector. Now define 
    \begin{equation*}
        \theta (z) \,:=\, z^{-M \lambda} \,=\, \big(z^{-M \cdot \lambda_1}, \dots, z^{-M \cdot \lambda_N}\big),
    \end{equation*}
    with $M := \max_{i \in [1,N]} \mathrm{deg} \,w_i$ being the
    maximal degree occurring among the polynomials $w_i$.

    By definition, we have $\pi = \pi \circ \theta(z)$ for
    every $z \in \mathbb H$, considering $\theta(z)$ as the map $\mathbb
    C^N \to \mathbb C^N$, $\xi \mapsto \theta(z) \cdot \xi$. It follows
    that
    \begin{equation*}
        D\pi\vert_{w(z)} \,=\, D\pi \vert_{\theta(z) \cdot w(z)} \circ D
        \theta(z)\vert_{w(z)},
    \end{equation*}
    and hence
    \begin{equation*}
        D\pi\vert_{w(z)} \tilde \eta_1(z) \,=\, (D\pi \vert_{\theta(z)
            \cdot w(z)} \circ D \theta(z)\vert_{w(z)})\tilde \eta_1(z).
    \end{equation*}
    Now $D\theta(z) \vert_{w(z)}$ converges to $0$ for $z \to
    \infty$ by the decay property of $\tau$. Moreover, since both $z \mapsto
    D\pi\vert_{\theta(z) \cdot w(z)}$ and $z \mapsto \tilde \eta_1 (z)$
    are bounded (the first because $\theta \cdot w$ only takes values in
    some compact subset of $\mathbb C^N$ by the decay of $\theta$, and the
    second because it converges), we get $D\pi_{w(z)} \tilde \eta(z)
    \to 0$.
\end{proof}

This completes the proof of the claim and hence of the disc part of Theorem
\ref{thm::transversality}.

\subsection{Regularity for spheres.} The proof of the sphere part of
Theorem \ref{thm::transversality} follows the same pattern as that of
the disc part. Let $v: \mathbb C P^1 \to X$ be an arbitrary
holomorphic sphere. As in the disc case, proving regularity of $J_0$
at $v$ amounts to showing that the Dolbeault operator
$\overline \partial: C^\infty(v^*TX) \to \Omega^{0,1}(v^*TX) $ has
$\mathrm{coker} \,\overline
\partial = 0$. By Proposition \ref{prop::fine-resolution}, which also
holds for closed Riemann surfaces, this translates into $H^1(\mathbb C
P^1; \mathcal O(v^*TX)) = 0$, where $\mathcal O(v^*TX)$ denotes the
sheaf of holomorphic sections of $v^*TX$.

By arguments similar to the ones in the proof of Proposition
\ref{thm::hom_coords}, the restriction of $v$ to any subset $\mathbb C
P^1 \smallsetminus \{\mathrm{pt}\}$ admits a lift $w =
(w_1,\dots,w_N)$ to $\mathbb C^N$ with polynomial entries $w_i$.
One obtains the exact sequence of sheaves
\begin{equation*}
    0 \to \mathcal O(w^*(TOrb_{K_{\mathbb C}})) \to \mathcal O(w^*\mathbb C^N)
    \to \mathcal O(v^*TX) \to 0,
\end{equation*}
and then deduces from the resulting long exact sequence in cohomology
that the cover $\mathfrak V$ of $\mathbb C P^1$ given by $V_0 =
\mathbb C P^1 \smallsetminus \{\infty\} \quad \text{and} \quad V_1 =
\mathbb C P^1 \smallsetminus \{0\}$ is acyclic for $\mathcal O
(v^*TX)$. The vanishing of $H^1(\mathbb C P^1;\mathcal O(v^*TX)) \cong
\check{H}^1(\mathfrak V;\mathcal O(v^*TX)$ is then shown by the same
\v{C}ech type argument as in the previous section: Given a holomorphic
section $\sigma$ of $v^*TX$ over $V_0 \cap V_1$, lift it to a
holomorphic section of the trivial bundle, decompose it using its
Laurent expansion as $\sigma = \eta_0 + \eta_1$ such that $\eta_0$
extends over $0$ and $\eta_1$ extends over $\infty$, and then project
back to $v^*TX$.

\section{The Involution on the Moduli Space of Real Discs}
\label{sec:invol-moduli-spac}
Denote by $H_2^D(X,R)$ the image of the Hurewicz homomorphism
$\pi_2(X,R) \to H_2(X;\mathbb Z)$ and by $\sim$ the equivalence
relation on $H_2^D(X,R)$ given by identifying classes of the same
Maslov number:
\begin{equation*}
    \lambda \,\sim\, \lambda' \quad : \Longleftrightarrow \quad
    \mu(\lambda) = \mu(\lambda').
\end{equation*}
Let $\widetilde{\mathcal M}(\lambda)$ be the moduli space of
para\-metrised holomorphic discs $u: (D^2,\partial D^2) \to (X,R)$
representing the class $\lambda \in H_2^D(X,R)/\mathord\sim$. Then
define the corresponding moduli space of unparametrised discs with two
marked boundary points by
\begin{equation*}
    \mathcal M_2(\lambda) \,:=\, \widetilde{\mathcal M}(\lambda)
    / \mathrm{Aut}_{\pm 1}(D^2),
\end{equation*}
where $\mathrm{Aut}_{\pm}(D^2)$ denotes the subgroup of
$\mathrm{Aut}(D^2)$ consisting of those automorphism which fix $\pm 1
\in \partial D^2$. Furthermore, set
\begin{equation*}
    \mathcal M_2 \,:=\, \bigsqcup_{\lambda} \,\mathcal M_2(\lambda)
\end{equation*}
We commonly denote elements of $\mathcal M_2$ by the same name as
their representatives, which should cause no confusion. Note that it
makes sense to evaluate $u \in \mathcal M_2$ at the points $\pm 1$.

\subsection{The involution $\tau_*$ and its fixed points.}
\label{sec:involution}
In this section we will study the involution
\begin{eqnarray*}
    \tau_*: \mathcal M_2 &\to& \mathcal M_2,\\
    u &\mapsto& \tau \circ u \circ c,
\end{eqnarray*}
where $c: D^2 \to D^2$ denotes complex conjugation. Its fixed point set
was characterized by Fukaya--Oh--Ohta--Ono in \cite{FOOO--Chapter_8}
(which is part of the big FOOO preprint, but is not contained in the
book version) in terms of another map
\begin{equation*}
    \mathfrak D: \mathcal M_2 \to  \mathcal M_2.
\end{equation*}
We briefly describe its definition. Denote by $D_\pm^2 := \{z \in D^2 ~|~ \pm
\mathrm{Im}\,z \geq 0\}$ the upper respectively lower half-disc, and
let $\rho_\pm: D_\pm^2 \to D^2$ be conformal isomorphisms satisfying
$\rho_-(z) = \overline{\rho_+(\overline z)}$. (In the upper half-plane
model, in which $D_\pm^2$ correspond to $\mathbb H_\pm = \big\{ z \in
\mathbb H ~\big\vert~ \pm \mathrm{Re}\,z \geq 0 \big\}$, such maps are
given by $\mathbb H_\pm \to \mathbb H$, $z \mapsto \pm z^2$.) Given a
real disc $u \in \widetilde{\mathcal M}_2(\lambda)$, define a new real
disc $u' \in \widetilde{\mathcal M}_2(\lambda - \tau_*\lambda) =
\widetilde{\mathcal M}_2(2 \lambda)$ by
\begin{equation*}
    \check u \,:= \,
    \begin{cases}
        \,u(\rho_-(z)),& z \in D_-^2\\
        \,\tau(u(\overline{\rho_+(z))}),& z \in D_+^2.
    \end{cases}
\end{equation*}
It is shown in \cite{FOOO--Chapter_8} (Lemma 40.4) that $\mathfrak D$
descends to a map $\mathfrak D: \mathcal M_2 \to \mathcal M_2$.

The fixed point set of $\tau_*$ is simply the image of $\mathfrak D$:

\begin{lem}[see Lemma 40.6 in \cite{FOOO--Chapter_8}] 
    \label{lem::charact-of-fixed-points}
    Given any $u \in \mathcal M_2$ which is fixed by $\tau_*$, there
    exists some $u' \in \mathcal M_2$ such that $\mathfrak D u' =
    u$.
\end{lem}

\begin{figure}[tbp]
    \centering
    \includegraphics[scale=0.3]{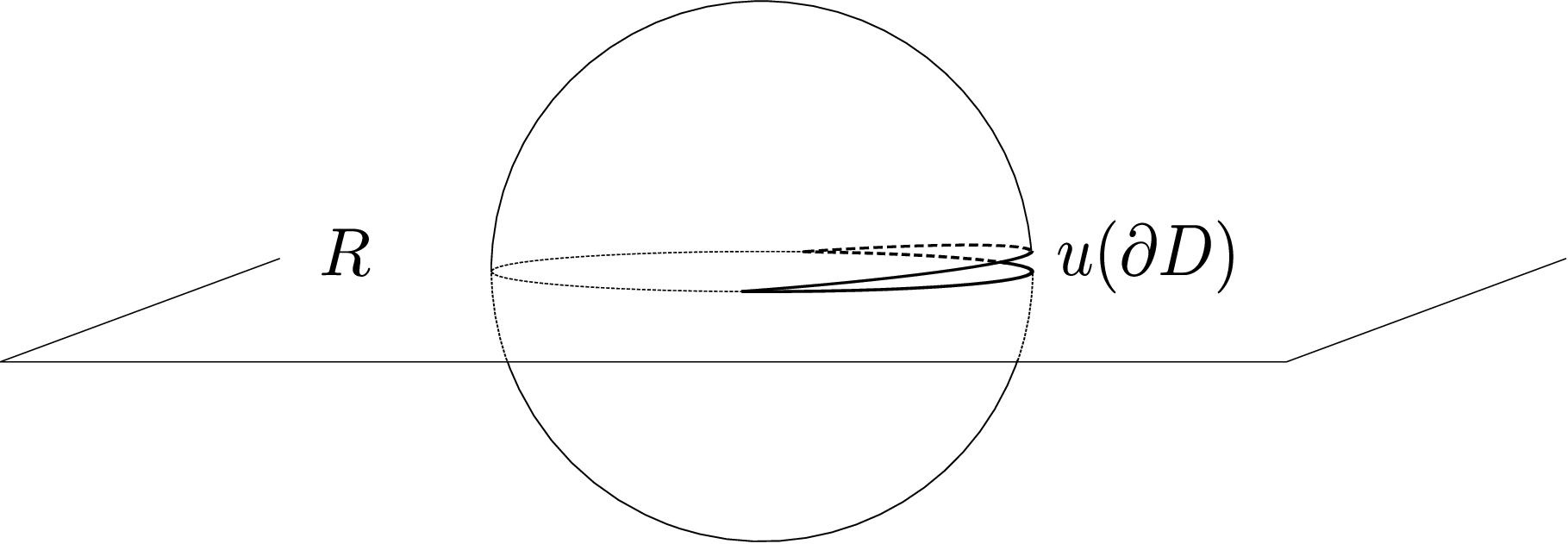}
    \caption{A Pac-Man disc}
    \label{fig:Pac-Man_disc}
\end{figure}

The fixed points can be pictured as illustrated in Figure
\ref{fig:Pac-Man_disc} (the mouth of the Pac-Man, which corresponds to
$\partial D^2$, should really be entirely closed).

\subsection{A closer look at the fixed points}
\label{sec::fixed-points-2}
The restrictions of a real rational curve $v: \mathbb C P^1 \to X$ to
the upper respectively lower hemispheres yield two real discs
\begin{equation*}
    v_\pm:  (D^2,\partial D^2) \to (X,R),
\end{equation*}
which we call the \emph{upper and lower halves} of $v$. In general the
$v_\pm$ are not distinct as elements of $\mathcal M_2$.

Let now $u \in \mathcal M_2$ be a fixed point of $\tau_*$, and let $u'
\in \mathcal M_2$ be such that $\mathfrak Du' = u$. Since every real
disc can be obtained as one half of its double, we have in particular
$u' = v_\pm$ for a real rational curve $v$ and can write
\begin{equation*}
    u = \mathfrak D v_\pm.
\end{equation*}
The next proposition formulates a dichotomy based on whether $v$ is
\emph{simple} or not. Recall that simplicity for a holomorphic sphere
$v: \mathbb C P^1 \to X$ means that it is not multiply covered,
i.e. it cannot be written as $v \,=\, v' \circ \phi$, where $v':
\mathbb C P^1 \to X$ is another holomorphic sphere and $\phi: \mathbb
C P^1 \to \mathbb C P^1$ is a holomorphic branched cover of degree $>
1$.

\begin{prop}
    \label{prop::dichotomy-for-fixed-points}
    Let $v: \mathbb C P^1 \to X$ be a real rational curve.
    \begin{enumerate}
        \renewcommand{\theenumi}{\roman{enumi}}
    \item If $v$ is simple, then $\mathfrak Dv_+$ and $\mathfrak
        Dv_-$ are distinct as elements of $\mathcal M_2$.
    \item If $v$ is non-simple, then there exist real discs $u_\pm$
        satisfying
        \begin{itemize}
        \item $\mu (u_\pm) < \mu(\mathfrak D v_\pm)$,
        \item $(u_\pm)(D^2) = (\mathfrak Dv_\pm)(D^2)$.
        \end{itemize}

    \end{enumerate}
\end{prop}

\begin{proof}
    (i) Suppose that $v$ is simple. In the usual identification
    $\mathbb C P^1 \smallsetminus \{\infty\} \,\cong\, \mathbb C$, the
    discs $v_\pm$ are simply the restrictions of $v$ to the upper and
    lower half-planes, and the $\mathfrak D v_\pm =: \check v_{\pm}$
    are given by $z \mapsto v(\pm z^2)$ for $z \in \mathbb H \subset
    D^2$. The images of $\partial D^2$ under $\check v_\pm$ are hence
    \begin{equation*}
        \check v_+ (\partial D^2) \,=\, v \big(\overline {\mathbb
            R}_+\big) \quad \text{and} \quad \check v_-
        (\partial D^2) \,=\, v \big(\overline {\mathbb R}_-\big),
    \end{equation*}
    where $\overline{\mathbb R}_{\pm} = \mathbb R_\pm \cup \{\infty\}$
    denote the extended positive and negative half-lines.

    Suppose that $\check v_+$ and $\check v_-$ define the same element
    of $\mathcal M_2$. Then it follows in particular that $\check
    v_+(\partial D^2) = \check v_-(\partial D^2)$, and hence
    $v(\overline{\mathbb R}_+) = v(\overline{\mathbb R}_-)$. So the
    set of non-injective points of $v$ contains at least all of
    $\mathbb R P^1$, contradicting the fact that this set is at most
    countable for any simple curve, see Proposition 2.5.1 in
    \cite{mcduff-salamon04}.

    (ii) If $v$ is non-simple, we can by definition write $v = v'
    \circ \phi$ with $v': \mathbb C P^1 \to X$ another holomorphic
    sphere and a $\phi: \mathbb C P^1 \to \mathbb C P^1$ a branched
    cover of degree $> 1$.
    \begin{claim*}
        Such $v'$ can be chosen among real rational curves.
    \end{claim*}
    Observe that the statement of part (ii) of the proposition follows
    from the claim, because we can then define the required discs
    $u_\pm$ by
    \begin{equation*}
        u_\pm := v_\pm',
    \end{equation*}
    i.e. as the upper respectively lower halves of $v'$. These clearly
    have the required properties.
    
    \emph{Proof of the Claim.} We first construct a Riemann surface
    $(\Sigma,j)$ from the image of $v$, together with a
    $j$--anti-holomorphic involution $\nu: \Sigma \to \Sigma$ and a
    holomorphic map $f: \Sigma \to X$ satisfying $f = \tau \circ f
    \circ \nu$. Then we argue that $\Sigma$ is biholomorphic to
    $\mathbb C P^1$ by a biholomorphism which is intertwines $\nu$ and
    complex conjugation.

    As for the construction of $\Sigma$, we follow closely the
    exposition in the proof of Proposition 2.5.1 in
    \cite{mcduff-salamon04} (and refer there for more details), where
    essentially the same appears. Denote by $C \subset \mathrm{Im}(v)$
    the set of critical values of $v$, and let $B$ be the subset of
    $\mathrm{Im}(v) \smallsetminus C$ consisting of all those points
    where distinct branches of $v$ meet. That is, a point $x \in
    \mathrm{Im}(v)$ is in $B$ iff there exist distinct $z_0,z_1 \in
    \mathbb C P^1 \smallsetminus v^{-1}(C)$ such that $v(z_1) =
    v(z_2)$ and such that for all sufficiently small neighbourhoods
    $U_i$ of $z_i$ we have $v(U_0) \neq v(U_1)$. Since $B$ is a
    discrete subset of $\mathrm{Im}(v)$, it follows that
    \begin{equation*}
        \Sigma' \,:=\, \mathrm{Im}(v) \smallsetminus
        (B \cup C)
    \end{equation*}
    is an embedded submanifold. We denote the embedding by $\iota:
    \Sigma' \to X$. $\Sigma'$ carries a unique complex
    structure $j'$ such that $\iota$ is holomorphic. Furthermore,
    $\Sigma'$ has two types of ends, namely those corresponding to the
    points in $B$ (a finite number for each $x \in B$, one for
    each branch of $v$ through $x$), and those corresponding to the
    points in $C$. By adding one point to $\Sigma'$ for each end and
    extending the complex structure, we obtain a closed Riemann
    surface $(\Sigma,j)$, and the embedding $\iota$ extends to a
    holomorphic map $f: \Sigma \to X$.
    
    The Riemann surface $(\Sigma',j')$ carries an anti-holomorphic
    involution $\nu$ induced by $\tau$, i.e. given by $\nu =
    \iota^{-1}\circ \tau \circ \iota$, which extends to an
    anti-holomorphic involution of $(\Sigma,j)$. Furthermore, $f$
    satisfies $f = \tau \circ f \circ \sigma$.

    By construction, $\Sigma$ admits a branched cover $\mathbb C P^1
    \to \Sigma$ and must hence itself be biholomorphic to $\mathbb C
    P^1$. We argue that there exists a biholomorphism $\psi: \mathbb C
    P^1 \to \Sigma$ intertwining complex conjugation and $\nu$.
    Indeed, this follows from the fact that complex conjugation is the
    only anti-holomorphic involution on $\mathbb C P^1$ (up to
    conjugation with an automorphism) having a non-discrete set of
    fixed points, which is also the case for $\nu$ by construction.

    Finally, setting
    \begin{equation*}
        v' \,:=\, f \circ \psi : \mathbb C P^1 \to X
    \end{equation*}
    produces the required real rational curve, which ends the proof of
    the claim and hence of the proposition.
\end{proof}

\section{The Pearly Moduli Spaces}
\label{sec:transv-eval-maps}
In this section we recall from
\cite{biran-cornea_ALagQH09,biran-cornea_QuantumStrForLag07} the
definitions of the moduli spaces appearing in the construction of the
quantum homology rings $QH(R;\Lambda_R)$ and $QH(X;\Lambda_X)$, and
show that the necessary auxiliary structures may be chosen among
$\tau$--(anti)-invariant ones, which is essential for the proof of
Theorem \ref{Thm::Wideness_of_R} to work.

\subsection{Definitions}
\label{sec:definitions-pearly-moduli-spaces}
\begin{figure}[tbp]
    \centering
    \includegraphics[scale=0.45]{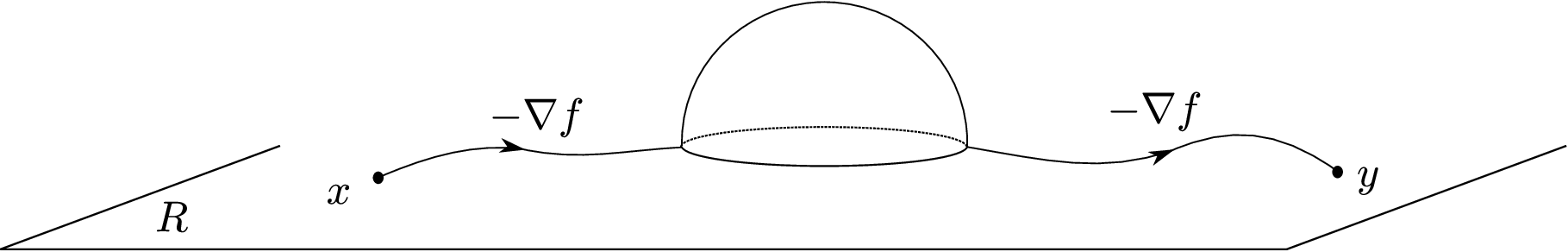}\\~\\~\\
    \includegraphics[scale=0.45]{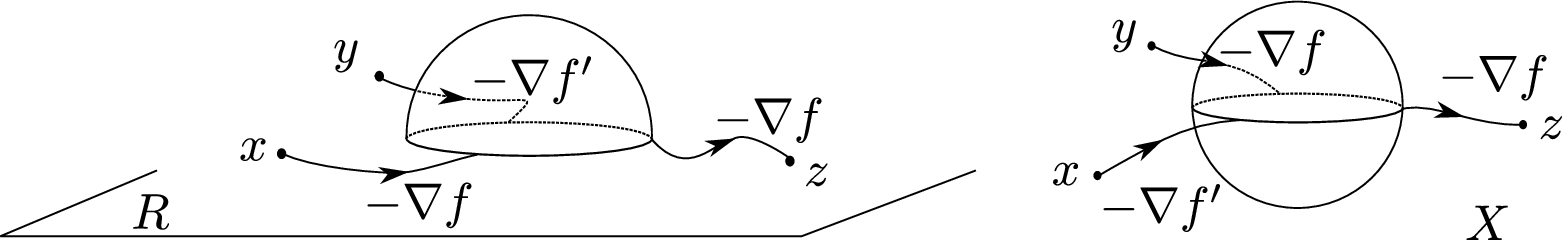}
    \caption{Elements of the pearly moduli spaces}
    \label{fig:pearly-traj}
\end{figure}

Fix a triple $\mathcal D = (f,g,J)$, where $f: X \to \mathbb R $ is a
Morse function, $g$ a Riemannian metric on $X$ and $J$ an almost
complex structure on $X$. Denote for now by $f$ and $g$ also the
restrictions to $R$, and by $\Phi_t$ the negative gradient flow of
$f$. We consider three different types of moduli spaces, which appear
in the definitions of the Lagrangian quantum differential and the quantum
products on $QH(R;\Lambda_R)$ respectively
$QH(X;\Lambda_X)$. (A graphical depiction of the elements of these
moduli spaces is provided in Figure \ref{fig:pearly-traj}.)

\subsubsection{Pearly trajectories} Given points $x,y \in R$ (not
necessarily critical) and a class $\lambda \in H_2^D(X,R)
/\mathord\sim$ (recall that $\sim$ identifies classes of the same
Maslov number), denote by $ \mathcal P(x,y,\lambda;\mathcal D) $ the
set of all tuples $(u_1,\dots,u_\ell) \in \mathcal{M}_2^{\times \ell}$
such that
\begin{enumerate}
\item All $u_i$ are non-constant,
\item $\Phi_{t_0}(u_1(-1)) = x$ for some $-\infty \leq t_0 < 0$,
\item $\Phi_{t_i}(u_i(+1)) = u_{i+1}(-1)$ for some
    $0<t_i<\infty$, $i \in [1,\ell-1]$,
\item $\Phi_{t_\ell}(u_\ell(+1)) = y$ for some $0 < t_{\ell} \leq \infty$,
\item $[u_1] + \dots + [u_\ell] = \lambda$.
\end{enumerate}
Here $\ell$ is not supposed to be fixed. If $x$ and $y$ are critical
points of $f$, the virtual dimension of this space is
\begin{equation*}
    \dim_{\mathrm{virt}} \mathcal P(x,y,\lambda;\mathcal D) = \vert x \vert - \vert y \vert    + \mu(\lambda) - 1,
\end{equation*}
where $\vert \cdot \vert$ denotes the Morse index of $f: R \to \mathbb
R$.

\subsubsection{Discy Y--pearls} Choose a second Morse function $f': X
\to \mathbb R $ such that all the (un)stable submanifolds of $f$ and
$f'$ intersect transversely. Let $x,z \in \mathrm{Crit}\,f$, $y \in
\mathrm{Crit}\,f'$, $\lambda \in H_2^D(M,L) /\mathord\sim$. Then
define $ \mathcal Y^D(x,y,z,\lambda;\mathcal D, \mathcal D')$
to be the set of all tuples $(\mathbf{u}, \mathbf{u}',\mathbf{u}'',u_c)$
where
\begin{enumerate}
\item $u_c: (D,\partial D) \to (X,R)$ is a real disc, possibly constant,

\item $\mathbf{u} \in \mathcal P(x, u_c(e^{-2 \pi i/3}),\eta;f,g,J)$,

\item $\mathbf{u}' \in \mathcal P(y, u_c(e^{2 \pi
        i/3}),\eta';f',g,J)$,
\item $\mathbf{u}'' \in \mathcal P(u_c(1),z,\eta'';f,g,J)$,
\item $\eta,\eta',\eta'' \in H_2^D(x,R)/\mathord\sim$ satisfy
    $\eta + \eta' + \eta'' + [u_c] = \lambda$.
\end{enumerate}
We will later refer to $u_c$ as the \emph{central disc}. The virtual
dimension in this case is
\begin{equation*}
    \dim_{\mathrm{virt}} \mathcal Y^D(x,y,z,\lambda;\mathcal D,
    \mathcal D') = \vert x \vert + \vert y \vert -
    \vert z \vert + \mu(\lambda) - n
\end{equation*}
with $\vert \cdot \vert$ denoting again the Morse index of $f: R \to
\mathbb R$.

\subsubsection{Spherical Y--pearls.}
Denote by $H_2^S(X)$ the image of the Hurewicz homomorphism $\pi_2(X)
\to H_2(X;\mathbb Z)$ and by $\sim$ the equivalence relation given by
identifying classes of the same Chern number. Again, let $f': X \to
\mathbb R$ be a second Morse function and denote its negative gradient
flow by $\Phi'_t$. Given some $x,z \in \mathrm{Crit}\,f$ and $y \in
\mathrm{Crit}\,f'$ and $\lambda \in H_2^S(X) / \mathord \sim$, let
$\mathcal Y^S (x,y,z,\lambda;\mathcal D, \mathcal D')$ be the space of
all holomorphic spheres $v: \mathbb C P^1 \to X$ such that
\begin{enumerate}
\item $\Phi_{-\infty} (v(e^{-2 \pi i/3})) = x$,
\item $\Phi'_{-\infty} (v(e^{2 \pi i/3})) = y$,
\item $\Phi_{+\infty}(v(1)) = z$,
\item $[v] = \lambda$.
\end{enumerate}
The virtual dimension of this space is
\begin{equation*}
    \dim_{\mathrm{virt}} \mathcal Y^S (x,y,z,\lambda;\mathcal D,
    \mathcal D') = \vert x \vert + \vert y \vert - \vert z \vert +
    2c_1(\lambda) -2n,
\end{equation*}
where now $\vert \cdot \vert$ denotes the Morse index of $f,f': X \to
\mathbb R $.

\subsection{Regularity}
\label{sec:transversality-eval-maps}
This subsection serves to justify working with the following data: The
standard $\tau$--anti-invariant complex structure $J_0$ on $X$ and a
generic pair $(f_\xi, g)$, where $f_\xi: X \to \mathbb R$ is defined
by $f_\xi(x) = \langle \mu, \xi \rangle$ with some $\xi \in \mathfrak
q$ as before, and $g \in \mathcal G_\tau$ is a $\tau$--invariant Riemannian metric on
$X$. More precisely, we prove the following proposition.

\begin{prop}
    \label{prop:regularity-mod-spaces}
    There exist second category subsets of $\mathfrak q \times
    \mathcal G_\tau$ resp. $\mathfrak q \times \mathfrak q \times
    \mathcal G_\tau$ such that for all $(\xi,g)$ resp. $(\xi,\xi',g)$
    in these subsets all moduli spaces $\mathcal
    P(x,y,\lambda;\mathcal D)$ resp. $\mathcal Y^D(x,y,z,\lambda;\mathcal
    D, \mathcal D')$ and $\mathcal Y^S(x,y,z,\lambda;\mathcal D,
    \mathcal D')$ defined using data triples $\mathcal D =
    (f_\xi,g,J_0)$ resp. pairs of data triples $\mathcal D=
    (f_\xi,g,J_0)$ and $\mathcal D' = (\tilde f_{\xi'},g,J_0)$ are
    smooth manifolds of the respective virtual dimension.\footnote{As
        in the case of the Morse intersection product (see the end of
        Section \ref{sec:homology-of-R}), we use as the Morse function
        in the second triple used for the definition of the quantum
        product a slightly perturbed version $\tilde f_{\xi'}$ of
        $f_{\xi'}$, such as to ensure that the critical sets of
        $f_\xi$ and $\tilde f_{\xi'}$ are disjoint.}
\end{prop}

\begin{proof}
    The proof uses only standard methods. We give an outline for the
    case $\mathcal P(x,y,z,\lambda;\mathcal D)$, the other cases being
    similar. Let $\boldsymbol{\lambda} =
    (\lambda_1,\dots,\lambda_\ell)$ be a vector of non-zero homology
    classes $\lambda_i \in H_2^D(X,R)$ and set
    \begin{equation*}
        \mathcal {M}_2(\boldsymbol\lambda) := \mathcal M_2(\lambda_1) \times \dots
        \times \mathcal M_2(\lambda_\ell),
    \end{equation*}
    which is a finite dimensional manifold. Then set
    \begin{equation*}
        \mathcal W^{1,2}_\ell := W^{1,2}((-\infty,0],R) \times
        (W^{1,2}([0,1],R))^{\times (\ell-1)} \times
        W^{1,2}([0,\infty),R),
    \end{equation*}
    where $W^{1,2}(I,R)$ denotes the space of paths $\gamma: I \to R$
    of Sobolev class $W^{1,2}$. It is well-known that $W^{1,2}(I,R)$
    is a smooth Banach manifold, and hence $\mathcal W^{1,2}_\ell$ is
    so, too. Consider now the map
    \begin{equation*}
        EV: \mathcal M_2(\boldsymbol \lambda) \times \mathcal W^{1,2}_\ell
        \, \to \, \mathcal M_2(\boldsymbol \lambda) \times R^{2\ell +2},
    \end{equation*}
    \begin{equation*}
        (\mathbf{u},\boldsymbol \gamma) \, \mapsto\, \Big(\mathbf{u}, (\gamma_0(-\infty),\gamma_0(0),\gamma_1(0),\dots,\gamma_\ell(\infty))\Big).
    \end{equation*}
    This map is a submersion, so that the preimage under $EV$ of every
    submanifold of $\mathcal M_2(\boldsymbol \lambda) \times
    R^{2\ell+2}$ is a submanifold. In particular, the set
    \begin{equation*}
        \mathcal X := EV^{-1}\Big\{\big(\mathbf{u},(x,u_1(-1),\dots,u_\ell(+1),y)\big) \in \mathcal M_2(\boldsymbol \lambda) \times R^{2\ell+2} ~\Big|~ u
        \in \mathcal{M}_2(\boldsymbol \lambda)\Big\}
    \end{equation*}
    is a submanifold of $\mathcal M_2(\boldsymbol \lambda) \times
    \mathcal W^{1,2}_\ell$, because the set on the right-hand side is
    clearly a submanifold of $\mathcal M_2(\boldsymbol \lambda) \times
    R^{2\ell+2}$.

    Fix now some $k > 0$ and let $\mathcal G_\tau^k$ be the set of
    $\tau$--invariant metrics of class $C^k$, which is a Banach
    manifold. Consider then the Banach bundle $\mathcal E \to \mathcal
    X \times \mathfrak q \times \mathcal G_\tau^k$ whose fibre at
    $(\mathbf{u},\boldsymbol \gamma, \xi,g)$ is
    \begin{equation*}
        \mathcal E_{(\mathbf{u},\boldsymbol
            \gamma, \xi,g)} = L^2 (\boldsymbol \gamma^*TR).
    \end{equation*}
    Define a section $\psi$ of this bundle by
    \begin{equation*}
        \psi (\mathbf{u},\boldsymbol
        \gamma, \xi,g ) \,:=\, \boldsymbol \gamma' - V_\xi \circ
        \boldsymbol \gamma \,:=\, (\gamma_0' -
        V_\xi \circ \gamma_0,\dots,\gamma_\ell' - V_\xi
        \circ \gamma_\ell),
    \end{equation*}
    where $V_\xi$ denotes the negative gradient vector field of
    $f_\xi$ with respect to $g$. Observe that for given $(\xi,g)$, there
    is a natural 1--1 correspondence
    \begin{equation*}
        \big\{(\mathbf{u},\boldsymbol \gamma)\in \mathcal X ~\big|~ \psi(\mathbf{u},\boldsymbol
        \gamma,\xi,g) = 0\big\} \, \cong \, \mathcal P(x,y,\boldsymbol
        \lambda; (f_\xi,g,J_0)),
    \end{equation*}
    where the set on the right-hand side denotes the space of pearly
    trajectories whose discs represent the vector $\boldsymbol
    \lambda$.
    
    \begin{claim*} For every $(\mathbf{u}, \boldsymbol \gamma, \xi, g)
        \in \psi^{-1}(0)$, the following assertions hold true:
        \renewcommand{\theenumi}{\roman{enumi}}
        \begin{enumerate}
        \item The vertical differential $D^v \psi :T
            _{(\mathbf{u},\boldsymbol\gamma,\xi,g)}(\mathcal X \times
            \mathfrak q \times T_g \mathcal G_\tau^k) \to \mathcal
            E_{(\mathbf{u},\boldsymbol\gamma,\xi,g)}$ is onto.
        \item The restricted vertical differential $D_1^v \psi: T
            _{(\mathbf{u},\boldsymbol\gamma)}\mathcal X \to \mathcal
            E_{(\mathbf{u},\boldsymbol\gamma,\xi,g)}$ is a Fredholm
            operator of index $\vert x \vert - \vert y \vert +
            \mu(\boldsymbol \lambda) - 1$.
        \end{enumerate}
    \end{claim*}

    Once this is proven, the statement of the proposition will follow
    by first applying a standard theorem (see
    e.g. \cite{Hutchings--Notes_on_Morse_Homology_2002}), which says
    that under these conditions the set of $(\xi,g) \in \mathfrak q
    \times \mathcal G_\tau^k$ for which $\{(\mathbf{u}, \boldsymbol
    \gamma) ~|~ \psi(\mathbf{u},\boldsymbol \gamma,\xi,g) = 0 \}
    \subset \mathcal M_2(\boldsymbol \lambda) \times \mathcal
    W^{1,2}_\ell$ is a smooth submanifold of the right dimension is of
    second category, and then using a procedure due to Taubes to get
    from $C^k$ to $C^\infty$ (see the proof of Theorem 3.1.5 (ii) in
    \cite{mcduff-salamon04}).

    We omit the proof of assertion (ii), which is standard. As for
    (i), note that $D^v \psi$ evaluated on a tangent vector of the
    form $(0,0,0,\dot g)$ with $\dot g \in T_g \mathcal G_\tau^k$
    yields
    \begin{equation*}
        D^v \psi \,(0,0,0,\dot g) \,=\, - \dot V_\xi \circ
        \boldsymbol \gamma,
    \end{equation*}
    where $\dot V_\xi = -g^{-1} \dot g (V_\xi)$ is the variation of
    the gradient vector field induced by $\dot g$ (here we view $g$,
    $\dot g$ as sections of $\mathrm{Hom}(TX,T^*X)$). Let $W \in
    L^2(\boldsymbol \gamma^*TR)$ be any element orthogonal to the
    range of $D_1^v\psi$ (which is a closed subspace of $\mathcal
    E_{(\mathbf{u},\boldsymbol \gamma,\xi,g)}$, since $D_1^v \psi$ is
    Fredholm). Then in particular
    \begin{equation*}
        \int_{\boldsymbol \gamma} g( \dot V_\xi,W ) \, ds \,= \,0
    \end{equation*}
    for all $\dot g \in T_g \mathcal G_\tau^k$. Now given any point $p
    \in \mathrm{im}\,\boldsymbol \gamma \in R$, one can easily find
    some $\dot g(p) \in \mathrm{Hom}(T_pR , T_p^*R)$ such that $\dot
    V_{f_\xi} (p) = W(p)$ for the corresponding $\dot V_\xi$. Here we
    assume $W$ to be continuous, which is legal by density of
    continuous sections in $L^2$. Extending $\dot g(p)$ by using a
    $\tau$--invariant cut-off function supported near $p$ yields an
    element $\dot g \in T_g \mathcal G_\tau^k$, and the fact that it
    pairs to zero with $W$ allows to conclude that $W(p)= 0$. This
    argument works because $\boldsymbol \gamma$ consists of pieces of
    gradient flow lines, which cannot intersect each other. So $W
    \equiv 0$, and hence $D_1^v \psi$ is surjective.

\end{proof}

\subsection{Compactness}
\label{sec:compactness}
To be able to count the elements of the zero-dimensional pearly moduli
spaces, we need the following proposition. 
\begin{prop}
    Let $\mathcal D, \mathcal D'$ be generic data triples in the sense
    of Proposition \ref{prop:regularity-mod-spaces}. Then all moduli
    spaces $\mathcal P(x,y,\lambda;\mathcal D)$, $\mathcal
    Y^D(x,y,z,\lambda;\mathcal D, \mathcal D')$ and $\mathcal
    Y^S(x,y,z,\lambda;\mathcal D, \mathcal D')$ of dimension zero are
    compact, i.e. finite sets.
\end{prop}

The proof of this proposition is the point where the assumptions of
monotonicity and $C_X \geq 2$ are needed. Given that these are
satisfied, it requires no arguments which are special to our
situation, but follows directly from Biran--Cornea's general
theory. The idea for proving compactness in dimension zero is that if
there were some sequence without a converging subsequence, then a
subsequence would have to Gromov-converge to a pearly configuration in
some other space, e.g. by bubbling off of a sphere. From this one
could extract a configuration in a pearly moduli space of
\emph{negative} virtual dimension (in our example by forgetting the
sphere), contradicting the fact that such a space must be empty by
Proposition \ref{prop:regularity-mod-spaces}. For a precise version of
this argument, see \cite{biran-cornea_QuantumStrForLag07}.
\newpage

\section{Proof of Theorem \ref{Thm::Wideness_of_R}}
\label{sec:proofs-main-thms}

\subsection{Proof of Theorem \ref{Thm::Wideness_of_R} (\ref{Thm-part::wideness})}
\label{sec:proof-thmA}
We briefly recall the definition of the quantum homology
$QH(R;\Lambda_R)$, referring to
\cite{biran-cornea_ALagQH09,biran-cornea_QuantumStrForLag07} for
details. Given a data triple $\mathcal D = (f,g,J)$, define the graded
vector space
\begin{equation*}
    \mathcal C(\mathcal D; \Lambda_R) \,=\, \mathbb Z_2 \langle \mathrm{Crit}\,f
    \rangle \otimes \Lambda_R
\end{equation*}
whose grading is induced by the Morse index of $f$ and the grading of
$\Lambda_R$. Consider the $\Lambda_R$--linear homomorphism $d:
\mathcal C(\mathcal D;\Lambda_R) \to \mathcal C(\mathcal D;\Lambda_R)$ given by
\begin{equation*}
    d x \,=\, \sum_{y,\lambda} \big(\#_{\mathbb Z_2} \mathcal
    P(x,y,\lambda)\big) y t^{\mu(\lambda)/N_R},
\end{equation*}
for $x \in \mathrm{Crit}\,f$, where the sum runs over all $y \in
\mathrm{Crit}\,f$ and all $\lambda \in H_2^D(X,R)/\mathord\sim$ such
that $\vert x \vert - \vert y \vert - \mu(\lambda) - 1 = 0$. Assuming
a generic choice of $\mathcal D$, the count makes sense, because then
$\mathcal P(x,y,\lambda)$ is a compact 0--dimensional manifold, i.e. a
finite set of points. Moreover, $d$ is a differential,
i.e. $d^2=0$. The quantum homology of $R$ is then defined
as
\begin{equation*}
    QH(R;\Lambda_R) \,:=\,QH(\mathcal D;\Lambda_R)\,:=\,  H(\mathcal C(\mathcal D;\Lambda_R),d)
\end{equation*}
and is independent of the data up to isomorphism, in the sense that
given another data triple $\mathcal D'$ there exists a chain map
$\mathcal C( \mathcal D;\Lambda_R) \to \mathcal C(\mathcal
D';\Lambda_R)$ inducing a canonical isomorphism
\begin{equation*}
   \Psi_{\mathcal D',\mathcal D}: QH(\mathcal  D;\Lambda_R) \to QH(\mathcal D';\Lambda_R).
\end{equation*}

We now specialize to the case that $(f_\xi,g_0,J_0)$ is a triple as in
Section \ref{sec:transv-eval-maps}, with $f_\xi \in \mathcal F_0$
chosen such as to make all the corresponding pearly moduli spaces
manifolds of the expected dimension. We know that we compute the right
quantum invariants when working with such data, because we can relate
the corresponding moduli spaces to those for arbitrary generic data by
a compact cobordism. This follows from the regularity results in
Section \ref{sec:transv-eval-maps} and the general theory of
Biran--Cornea.

\begin{claim*}
    For this choice of data $\mathcal D = (f_\xi,g_0,J_0)$, the corresponding
    homomorphism
    \begin{equation*}
        d: \mathcal C(\mathcal D;\Lambda_R) \to \mathcal C(\mathcal D;\Lambda_R)
    \end{equation*}
    vanishes.
\end{claim*}

The wideness statement of Theorem \ref{Thm::Wideness_of_R} follows
immediately from the claim: The vanishing of $d$ implies firstly that
$d$ is a differential, and secondly that the complex $(\mathcal
C(\mathcal D;\Lambda_R), d = 0)$ is identical to the Morse complex associated to
$f_\xi$. Hence the map
\begin{equation}
    \label{eq:morphism-QH-HM}
    QH(\mathcal D;\Lambda_R) \to HM(f_\xi) \otimes \Lambda_R
\end{equation}
induced by mapping a quantum homology class $[x]_{Q}$ represented by some $x \in
\mathrm{Crit}\,f_\xi$ to the class $[x]_{M}$ represented by $x$ in
Morse homology gives rise to an isomorphism
\begin{equation*}
    QH(R;\Lambda_R) \cong H(R; \mathbb Z_2) \otimes \Lambda_R.
\end{equation*}
We will explain below why this isomorphism is canonical, in the sense
that the isomorphisms of type \eqref{eq:morphism-QH-HM} are compatible
with the comparison isomorphisms relating different choices of data.

\begin{proof}[Proof of the Claim]
    Let $x,y \in \mathrm{Crit}\,f_\xi$, let $\lambda \in
    H_2^D(X,R)/\mathord \sim$ be such that $\vert x \vert - \vert y
    \vert + \mu(\lambda) = 0$ and consider the
    involution
    \begin{equation*}
        \tau_*^{\mathcal P}: \mathcal P(x,y,\lambda) \to
        \mathcal P(x,y,\lambda),\quad
        (u_1,\dots,u_k) \mapsto (\tau_*u_1,\dots,\tau_*u_k),
    \end{equation*}
    on the space of pearly trajectories, where on the right $\tau_*$
    denotes the involution on $\mathcal M_2$ defined in Section
    \ref{sec:invol-moduli-spac}. Observe that
    \begin{equation*}
        \#_{\mathbb Z_2} \mathcal P(x,y,\lambda) =
        \#_{\mathbb Z_2} \mathrm{Fix}\,\tau_*^{\mathcal P}.
    \end{equation*}
    Hence $d = 0$ will follow if we can show that the number of fixed
    points of $\tau_*^{\mathcal P}$ is even. So suppose that
    $(u_1,\dots,u_k) \in \mathcal P(x,y,\lambda)$ is a fixed point of
    $\tau_*^{\mathcal P}$, which is equivalent to all $u_i$ being
    fixed points of $\tau_*: \mathcal M_2 \to \mathcal M_2$. By Lemma
    \ref{lem::charact-of-fixed-points} and the discussion in Section
    \ref{sec::fixed-points-2} we can write
    \begin{equation*}
        (u_1,\dots,u_k) := (\mathfrak Dv_{1,+}, \dots, \mathfrak Dv_{k,+})
    \end{equation*}
    with real rational curves $v_i: \mathbb C P^1 \to X$ (recall that
    $v_+$ denotes the upper half of $v: \mathbb C P^1 \to X$). We will
    argue below that all $v_i$ have to be simple. Assuming this for
    the moment, consider the pearly trajectory
    \begin{equation*}
        (u_1',\dots,u_k') := (\mathfrak Dv_{1,-}, \dots, \mathfrak Dv_{k,-}).
    \end{equation*}
    By Proposition \ref{prop::dichotomy-for-fixed-points} (i), we have
    $\mathfrak Dv_{i,+} \neq \mathfrak Dv_{i,-}$ in $\mathcal M_2$ and
    consequently
    \begin{equation*}
        (u_1,\dots,u_k) \,\neq\, (u_1',\dots,u_k')
    \end{equation*}
    in $\mathcal P(x,y,\lambda)$. This shows that all fixed
    points of $\tau_*^{\mathcal P}$ occur in pairs, so that their
    total number is even.

    It remains to justify the simplicity assumption. Consider an
    element of $\mathcal P(x,y,\lambda)$ containing a real disc $u =
    \mathfrak D v_+$ fixed by $\tau_*:\mathcal M_2 \to \mathcal M_2$.
    If $v$ were non-simple, we could infer by Proposition
    \ref{prop::dichotomy-for-fixed-points} (ii) that there exists some
    real disc $u'$ with
    \begin{equation*}
        \mu(u') < \mu(u)\quad \text{and} \quad
        u' (D) = u(D).
    \end{equation*}
    Replacing $u$ by $u'$ would hence yield a new pearly trajectory
    representing a homology class $\lambda'$ of Maslov number strictly
    lower than $\mu(\lambda)$. But this is not possible, because under
    our genericity assumptions $\mathcal P(x,y,\lambda')$ is a
    manifold of dimension
    \begin{equation*}
        \vert x \vert - \vert y \vert + \mu(\lambda') - 1 < 0
    \end{equation*}
    and thus empty. We conclude that $v$ is simple.
\end{proof}

It remains to address the question of why the isomorphism $QH(R;\Lambda_R) \to
H(R;\mathbb Z_2) \otimes \Lambda_R$ to which the family of
isomorphisms \eqref{eq:morphism-QH-HM} gives rise is
\emph{canonical}. What we mean by that is that all diagrams of the
form
\begin{diagram}
    QH(\mathcal D;\Lambda_R) & \rTo{} &  HM(f) \otimes \Lambda_R \\
    \dTo{\Psi_{\mathcal D,\mathcal D'}} & &\dTo>{\Psi^M_{f,f'}}\\
    QH(\mathcal D';\Lambda_R)& \rTo{} & HM(f') \otimes \Lambda_R \\
\end{diagram}
commute, where the horizontal arrows are as in
\eqref{eq:morphism-QH-HM} and the vertical arrows are the comparison
isomorphisms in quantum respectively Morse homologies. In general
$\Psi_{\mathcal D',\mathcal D}: QH(\mathcal D;\Lambda_R) \to QH(\mathcal D';\Lambda_R)$ is
of the form
\begin{equation*}
    \Psi_{\mathcal D',
        \mathcal D} = \Psi_{f',f}^M + \Psi_{\mathcal D',\mathcal D}^Q t
\end{equation*}
where $\Psi_{\mathcal D',\mathcal D}^Q$
counts certain configurations containing holomorphic discs (see
\cite{biran-cornea_QuantumStrForLag07,biran-cornea_ALagQH09} for
details). So the potential source of non-canonicality is that
$\Psi_{\mathcal D',\mathcal D}^Q$ might not vanish. However, in our
specific case a cancellation argument as in the proof of the claim
shows
\begin{equation*}
    \Psi_{\mathcal D',\mathcal D}^Q = 0
\end{equation*}
and hence our isomorphism is canonical. This ends the proof of Theorem
\ref{Thm::Wideness_of_R} (\ref{Thm-part::wideness}). \hfill $\square$

\subsection{Proof of Theorem \ref{Thm::Wideness_of_R}
    (\ref{Thm-part::ring_str})}
\label{sec:proof-thmB}
As in the proof of the wideness part we work with data triples of the
form $\mathcal D = (f_\xi,g_0,J_0)$, and when the product is concerned
with a second triple $\mathcal D' = (\tilde f_{\xi'},g_0,J_0)$, where
$f_\xi, f_{\xi'} \in \mathcal F_0$. 

We already know that $QH(R;\Lambda_R)$ and $QH(X;\Lambda_X)$ are
isomorphic to $H(R;\mathbb Z_2)\otimes \Lambda_R$ respectively
$H(X;\mathbb Z_2)\otimes \Lambda_X$ and hence isomorphic to each other
\emph{as modules}. Explicitly, a (degree-doubling) isomorphism is
induced by
\begin{equation*}
    \mathcal C(R;\mathcal D) \to \mathcal
    C(X;\mathcal D), \quad x \mapsto x,
\end{equation*}
where $f_\xi: X \to \mathbb R $ is a generic Morse function of the
type described in Section \ref{sec:transv-eval-maps}, and by the
degree-doubling ring isomorphism
\begin{equation*}
    \Lambda_R \to \Lambda_X
\end{equation*}
determined by $t \mapsto q$ (recall that $\vert q \vert = -2C_L = -2N_R =
2 \vert t \vert$). Here we denote by $\mathcal C(R;\mathcal D) =
\big(\mathbb Z_2 \langle \mathrm{Crit}f_\xi \rangle \otimes \Lambda_R,
d=0\big)$ the Lagrangian quantum complex and by $\mathcal C(X;\mathcal
D)$ the complex $\big(\mathbb Z_2 \langle \mathrm{Crit}f_\xi \rangle
\otimes \Lambda_X, d=0\big)$ whose homology is $QH(X;\Lambda_X)$. It
follows from the discussion at the end of the proof of Theorem
\ref{Thm::Wideness_of_R} (\ref{Thm-part::wideness}) that this
isomorphism is canonical.

It remains to show that the isomorphism is compatible with the quantum
products on both sides. We recall from \cite{biran-cornea_ALagQH09,
    biran-cornea_QuantumStrForLag07} that the quantum product on
$QH(R;\Lambda_R)$ is by definition induced by the $\Lambda_{\mathbb
    Z_2}$--linear homomorphism
\begin{equation*}
    \mathcal C(R;\mathcal D) \otimes \mathcal
    C(R;\mathcal D') \to \mathcal C(R;\mathcal D), \quad x \otimes y \,\mapsto\, x \ast y,
\end{equation*}
given by
\begin{equation*}
    x \ast y \,:=\, \sum_{z,\lambda} \big( \#_{\mathbb Z_2} \mathcal
    Y^D(x,y,z,\lambda)\big) z t^{\mu(\lambda)/N_R},
\end{equation*}
where $x,y,z$ and $\lambda \in H_2^D(X,R)$ are such that $\vert x \vert +
\vert y \vert - \vert z \vert + \mu(\lambda) -n =0$. Similarly, the
quantum product on $QH(X;\Lambda_X)$ is defined by $\mathbb
Z_2$--counting elements in the moduli spaces
\begin{equation*}
    \mathcal Y^S(x,y,z;\lambda)
\end{equation*}
of spherical pearls described in Section \ref{sec:transv-eval-maps},
with critical points $x,y,z$ and $\lambda \in H_2^S(X)$.

Observe that there is a 1--1 correspondence between the sets $H_2^D(X,R)
/ \mathord\sim$ and $H_2^S(X)/\mathord\sim$ induced by mapping the
class $\lambda$ represented by some disc to the class $\lambda^\sharp$
represented by its double (recall that $\sim$ identifies classes of
the same Maslov respectively Chern numbers). The theorem will hence
follow if we can show
\begin{equation*}
    \#_{\mathbb Z_2} \mathcal Y^D(x,y,z,\lambda) = \#_{\mathbb Z_2} \mathcal
    Y^S(x,y,z,\lambda^\sharp).
\end{equation*}

To do so, we begin by decomposing the moduli space $\mathcal
Y^D(x,y,z,\lambda)$ as
\begin{equation*}
    \mathcal Y^D (x,y,z,\lambda) \,=\, \mathcal Y^D_0 (x,y,z,\lambda) \,\sqcup\, \mathcal Y^D_1(x,y,z,\lambda),
\end{equation*}
where $\mathcal Y^D_0(x,y,z,\lambda)$ is the subset consisting of all
those Y--pearls containing some non-constant non-central disc, so that
$\mathcal Y^D_1(x,y,z,\lambda)$ is the subset of those Y--pearls
containing at most a central disc.

\begin{claim}
    \label{claim1}
    $\#_{\mathbb Z_2} \mathcal Y^D_0 (x,y,z,\lambda)\, = \, 0$.
\end{claim}

\begin{proof}
    Consider the involution
    \begin{equation*}
        \tau_*^{\mathcal Y^S}: \mathcal Y^D(x,y,z,\lambda) \to \mathcal Y^D(x,y,z,\lambda)
    \end{equation*}
    given by reflecting all discs in a product pearl except the
    central one (note that reflecting the central one would lead to an
    element in the moduli space $\mathcal Y^S(y,x,z,\lambda)$ with $x$
    and $y$ reversed). Using precisely the same arguments as in the
    proof of Theorem \ref{Thm::Wideness_of_R} (i) one shows
    that all fixed points of $\tau_*^{\mathcal Y^S}$ occur in pairs,
    which proves the claim.
\end{proof}
\begin{claim}
   \label{claim2}
   $\#_{\mathbb Z_2} \mathcal Y^D_1(x,y,z,\lambda) = \#_{\mathbb Z_2}
   \mathcal Y^S(x,y,z,\lambda^{\sharp})$.
\end{claim}

\begin{proof}[Proof of Claim \ref{claim2}]
    There is an obvious map
    \begin{equation*}
        \mathcal Y^D_1(x,y,z,\lambda) \to
        \mathcal Y^S(x,y,z,\lambda^{\sharp})
    \end{equation*}
    given by replacing the central disc $u$ of any Y--pearl by its
    double $u^\sharp$.

    We argue that this map is injective. Suppose $u_0 \neq u_1$ and
    $u_0^\sharp = u_1^\sharp$. Then we must have $u_1 = \tau \circ u_0
    \circ c$ (the reflection of $u_0$), and in particular
    \begin{equation*}
        u_1(e^{-2\pi i /3}) \,=\, u_0(e^{2 \pi i/3}) \quad
        \text{and} \quad u_1(e^{2 \pi i/3}) \,=\,
        u_0(e^{-2 \pi i /3}).
    \end{equation*}
    Since $u_0$ appears in a Y--pearl, we may assume $u_0(e^{2 \pi
        i/3}) \neq u_0(e^{-2 \pi i/3})$ under our genericity assumptions,
    so that we obtain $u_0(e^{2\pi i/3}) \neq u_1(e^{2 \pi
        i /3})$ and hence $ u_0^\sharp(e^{2 \pi i/3}) \neq
    u_1^\sharp(e^{2 \pi i /3})$, contradicting $u_0^\sharp =
    u_1^\sharp$.

    The only Y--pearls which are not in the image of this map are
    those which are not $\tau$--invariant, in the sense that they
    contain a sphere which is not a real rational curve or a Morse
    trajectory which is not entirely contained in $R$. But for these
    it is easy to see that they occur in pairs and thus do not
    contribute to the quantum product.
\end{proof}

From the two claims we conclude that
\begin{eqnarray*}
    \#_{\mathbb Z_2} \mathcal Y^D(x,y,z;\lambda) &=& \#_{\mathbb Z_2}
    \mathcal Y^D_0(x,y,z;\lambda) + \#_{\mathbb Z_2} \mathcal
    Y^D_1(x,y,z;\lambda) \\
    &=& \#_{\mathbb Z_2} \mathcal Y^D_1(x,y,z;\lambda)\\
    &=& \#_{\mathbb Z_2} \mathcal Y^S(x,y,z;\lambda^\sharp),
\end{eqnarray*}
which ends the proof of Theorem \ref{Thm::Wideness_of_R} (\ref{Thm-part::ring_str}). \hfill $\square$

\bibliographystyle{hplain}
\bibliography{rl}

\end{document}